\documentclass{amsproc}

\usepackage{graphicx}
\usepackage{amsmath,amscd, amssymb}
\usepackage[all]{xy}

\newtheorem{theorem}{Theorem}[section]

\theoremstyle{definition}
\newtheorem{definition}[theorem]{Definition}

\theoremstyle{remark}

\numberwithin{equation}{section}

\newcommand{\integer}{\ensuremath{{\mathbb Z}}}

\newcommand{\real}{\ensuremath{{\mathbb R}}}
\newcommand{\complex}{\ensuremath{{\mathbb C}}}

\newcommand{\rational}{\ensuremath{{\mathbb Q}}}

\newcommand{\Torus}{\ensuremath{{\mathbb T}}}
\newcommand{\FF}{\ensuremath{{\mathcal F}}}
\newcommand{\proj}{\ensuremath{{\mathbb P}}}
\newcommand{\Torics}{\ensuremath{{\mathbf{Torics}}}}
\newcommand{\Stacks}{\ensuremath{{\mathbf{Stacks}}}}
\newcommand{\NCTop}{\ensuremath{{\mathbf{NCTop}}}}
\newcommand{\Man}{\ensuremath{{\mathbf{Man}}}}
\newcommand{\Fans}{\ensuremath{{\mathbf{Fans}}}}
\newcommand{\Cones}{\ensuremath{{\mathbf{Cones}}}}

\newcommand{\op}{\ensuremath{{\mathrm{op}}}}
\newcommand{\Polytopes}{\ensuremath{{\mathbf{Polytopes}}}}
\newcommand{\Groupoids}{\ensuremath{{\mathbf{Groupoids}}}}
\newcommand{\commAlg}{\ensuremath{{\mathbf{commAlg}}}}
\newcommand{\Diffeologies}{\ensuremath{{\mathbf{Diffeologies}}}}
\newcommand{\Top}{\ensuremath{{\mathbf{Top}}}}
\newcommand{\Alg}{\ensuremath{{\mathbf{Alg}}}}
\newcommand{\Mod}{\ensuremath{{\mathbf{Mod}}}}
\newcommand{\Sets}{\ensuremath{{\mathbf{Sets}}}}
\newcommand{\OO}{\ensuremath{{\mathcal O}}}
\newcommand{\HH}{\ensuremath{{\mathcal H}}}
\newcommand{\GG}{\ensuremath{{\mathcal G}}}
\newcommand{\XX}{\ensuremath{{\mathcal X}}}
\newcommand{\YY}{\ensuremath{{\mathcal Y}}}
\newcommand{\BB}{\ensuremath{{\mathcal B}}}
\newcommand{\TT}{\ensuremath{{\mathcal T}}}
\newcommand{\NN}{\ensuremath{{\mathcal N}}}
\newcommand{\MM}{\ensuremath{{\mathcal M}}}
\newcommand{\DD}{\ensuremath{{\mathcal D}}}
\newcommand{\SSS}{\ensuremath{{\mathcal S}}}
\newcommand{\redu}{\mathrm{red}}
\newcommand{\Img}{\mathrm{Im\ }}
\newcommand{\Ker}{\mathrm{Ker\ }}
\newcommand{\per}{\mathrm{per}}
\newcommand{\even}{\mathrm{even}}
\newcommand{\odd}{\mathrm{odd}}

\newcommand{\topo}{\mathrm{top}}

\newcommand{\Mat}{\mathrm{Mat}}
\newcommand{\Pro}{\ensuremath{{\mathbf{Pro}}}}
\newcommand{\Hom}{\mathrm{Hom}}
\newcommand{\dR}{\mathrm{dR}}
\newcommand{\Hol}{\mathrm{Hol}}
\newcommand{\Db}{\ensuremath{{\mathcal{D}^b}}}

\begin{document}
 
\title{Non-commutative toric varieties}

\author[L. Katzarkov]{Ludmil Katzarkov}
\address{Fakult\"{a}t f\"{u}r Mathematik, Nordbergstrasse 15, 1090 Wien, Austria}
\email{lkatzark@math.uci.edu}
\thanks{The first author was supported in part by NSF Grant \#000000.}

\author[E. Lupercio]{Ernesto Lupercio}
\address{Departamento de Matem\'{a}ticas, Cinvestav, Av. Instituto Polit\'{e}cnico Nacional \# 2508, Col. San Pedro Zacatenco, M\'{e}xico, D.F. CP 07360, M\'{e}xico.}
\email{lupercio@math.cinvestav.mx}
\thanks{The second author was partially supported by Conacyt.}

\author[L. Meersseman]{Laurent Meersseman}
\address{Centre de Recerca Matematica, Campus de Bellaterra, Edifici C - 08193 Bellaterra (Barcelona)}
\email{laurent.meersseman@u-bourgogne.fr}

\author[A. Verjovsky]{Alberto Verjovsky}
\address{Instituto de Matem\'{a}ticas, UNAM, Av. Universidad s/n. Col. Lomas de Chamilpa, C\'{o}digo Postal 62210, Cuernavaca, Morelos.}
\email{alberto@matcuer.unam.mx}
\thanks{The third author was partially supported by Conacyt.}

\date{March 6, 2013 and, in revised form, October 21, 2013.}


\begin{abstract}
In this short note, we introduce a new family of non-commutative spaces that we call \emph{non-commutative toric varieties}. We also briefly describe some of their main properties. The main technical tool in this investigation is a natural extension of LVM-theory for the irrational case. This note, aimed at graduate students, primarily consists  of a brief survey of the required background for the definition of a non-commutative toric variety, plus further results at the end. An definitve version os this paper is to appear in Contemporary Mathematics.
\end{abstract}

\maketitle

\section{Introduction}

Non-commutative toric varieties are to toric varieties what non-commutative tori are to tori and, as such, they can be interpreted in multiple ways: As (non-commutative) topological spaces, they are $C^*$-algebras associated to dense foliations, that is to say, deformations of the commutative $C^*$-algebras associated to tori. But while non-commutative tori correspond to foliations (deformations) on classical tori, non-commutative toric varieties correspond to foliations on so-called LVM-manifolds, which are certain intersections of real quadrics in complex projective spaces of a very explicit nature. The homotopy type of LVM-manifolds is described by moment angle complexes \cite{BP2002} that provides us with the basic topological information regarding non-commutative toric varieties. 

Non-commutative toric varieties can be interpreted as stratified parametrized families of non-commutative tori over a simple polytope: indeed, while for a classical toric variety the inverse image of a point $\mu^{-1}(p)$ under the moment map $\mu:X\to P$ is an ordinary torus, a non-commutative toric variety still has a moment map $\mu : \XX \to P$ (landing in an irrational simple polytope) but then $\mu^{-1}(p)$ is a non-commutative torus.

Additionally, non-commutative toric varieties have a complex analytic nature: they are non-commutative K\"ahler manifolds. To be able to express this behaviour, we have chosen the language of proto-complex manifolds (complex versions of diffeological spaces). This provides a functor-of-points approach to this family of non-commutative spaces helping us to form a moduli space of non-commutative toric varieties which is itself a complex orbifold.

The purpose of this note is to define rigurously non-commutative toric varieties. This note, primarily aimed at graduate students, consists of a brief survey of the required background for the definition of a non-commutative toric variety in sections 2,3,4,5, plus a very short announcement of the results of \cite{NCTorics}  at the end in sections 6 and 7.

In section 2 we review the classical theory of (commutative) toric varieties over the complex numbers; in section 3, we review non-commutative geometry \'{a} la Alaine Connes concentrating on the example of the non-commutative torus,  in section 4, we play with several variations on the concept of diffeology necessary for the definition;  in section 5, we review the most basic facts of LVM-theory which provides a definition of a classical toric variety as the leaf space of a certain holomorphic foliation on a LVM-manifold that generalizes the Calabi-Eckmann fibration. This constitutes the background part of the paper. 

We conclude with the results of \cite{NCTorics} in section 6, including a new proof of the McMullen conjecture for an arbitrary simplicial convex polytope; in section 7 we also offer further directions of research regarding the mirror symmetry program for non-commutative toric varieties.

We are very thankful to the referees of this paper for many comments that we believe have improved the exposition. We also would like to thank a conversation with Ulrike Tillmann that clarified some important points.

\section{Classical Toric Varieties}

\subsection{} In this section, we review the classical theory of toric varieties. The material of this section is well known, and it is presented here mainly to fix the notation. Some excellent references for this are: \cite{fulton1993introduction, audin2004torus, CoxToric2011, hamilton2007quantization}.

\subsection{} \emph{Compact toric varieties}\footnote{Of course, there are non-compact toric varieties, but in this paper, we will mostly ignore them.} $X$ can be defined as projective normal algebraic varieties over $\complex$ that  are compactifications of the complex algebraic torus $\Torus^n := (\complex^*)^n =\{ (t_1,\dots t_n) \colon t_i \neq 0 \}$, in such a manner that the natural action of the torus on itself extends to all of $X=\overline{\Torus^n}$ in an algebraic fashion. Notice that, by definition, the action of the real torus $T^n = (S^1)^n \subset \Torus^n$ extends to $X$. 

\subsection{} Toric varieties naturally carry the structure of K\"{a}hler manifolds; therefore, they can be studied as symplectic manifolds. All classical toric varieties that we consider in this first section are either smooth or orbifolds. 

It turns out that the action of the real torus $G=T^n$ on $X$ is Hamiltonian so we can consider the image of the moment map: $$\mu\colon X \to P \subset \real^n \cong \mathrm{Lie}(T^n).$$ The image of the moment map $P:=\mu(X)$ is a compact convex polytope whose combinatorial geometry contains all the information of the variety $X$. 

\subsection{} All compact convex polytopes $P$ that appear from the previous construction are characterized by the following property: \emph{at each vertex $v$ of $P$, there are exactly $n$ edges, and there exist shortest integer vectors $(v_1,\ldots,v_n)$ along each of the $n$ edges meeting at $v$, forming a $\rational$-basis of the subgroup $\rational^n \subset \real^n$}. The polytopes satisfying this rationality condition are called  \emph{Delzant polytopes}. The variety $X$ is smooth (and not only an orbifold) if and only if, additionally at each $v$, the vectors $(v_1,\ldots,v_n)$ form a $\integer$-basis of the lattice $\integer^n \subset \rational^n$, in which case $P$ is an \emph{integral Delzant polytope}.

\begin{figure}[htb]
\includegraphics[height=1in,width=1.3in,angle=0]{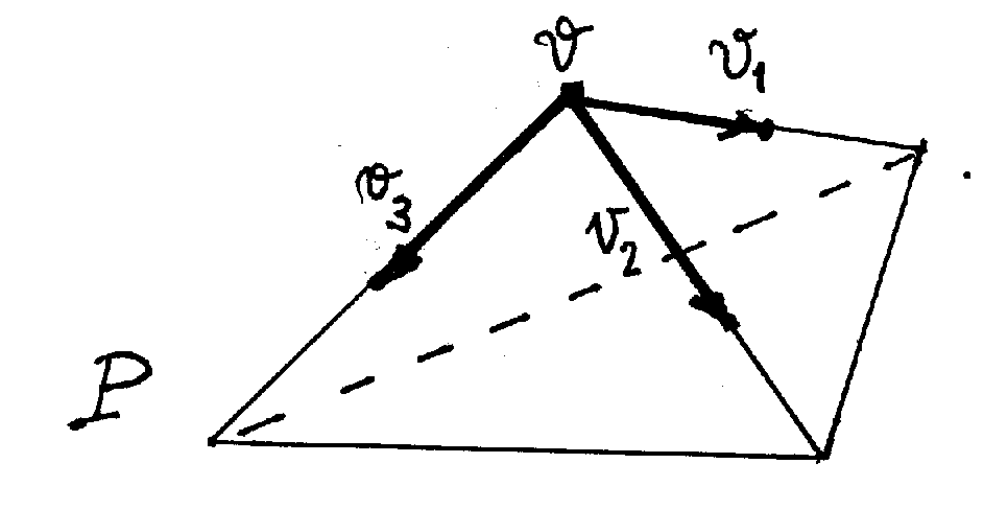}
\caption{A Delzant polytope.}
\label{firstfig}
\end{figure}

Notice that $\partial P$ provides a CW-decomposition of the sphere $S^n$ and that the boundary of the dual polytope $\partial P^\vee$ is combinatorially a triangulation of the sphere. The top dimensional faces $F_1,\ldots, F_N$ of codimension $1$ of $P$ are called its \emph{facets}, and we let $N$ be their number.

\subsection{} We can recover $X$ from $P$ using the following procedure: First, to $P\subset \real^n$, we assign a pair $(\FF, \rho)$ where $\FF$ encodes the combinatorics of $P$, and $\rho$ encodes all the geometry we need of the embedding $P \subset \real^n$:
\begin{itemize}
\item The family $\FF$ of subsets of $\{1,\ldots,N\}$ is defined as follows: $\emptyset \in \FF$;  otherwise $I\in\FF$ if and only if the intersection of facets  for indices in $I$ is non-empty, $\bigcap_{i \in I} F_i \neq \emptyset.$ Notice that $\FF$ is an abstract avatar for the triangulation of the sphere mentioned above.
\item The $N\times n$ integer matrix $\rho \colon \integer^N \to \integer^n$ assigns to the $i$-th vector $e_i$ of the canonical basis of $\integer^n$ the vector $\rho_i:=\rho(e_i)=(\rho_j^i)_{1\leq j \leq n} \in \real^n$ which is the shortest integer outward vector normal to the $i$-th facet of $P$. We will use the real matrix $\rho_\real := \rho \otimes \real \colon \real^N \to \real^n$, the complex matrix $\rho_\complex := \rho \otimes \complex \colon \complex^N \to \complex^n$,  as well as its (Lie theoretic) exponential $\hat{\rho} \colon \Torus^N \to \Torus^n$, $$\hat{\rho}(t_1,\ldots,t_N):=(t_1^{\rho^1_1} t_2^{\rho^2_1} \cdots t_N^{\rho^N_1}, \ldots, t_1^{\rho^1_n} t_2^{\rho^2_n} \cdots t_N^{\rho^N_n}).$$
\end{itemize}
Given this data $(\FF,\rho)$, we can assign an open set $U_\FF \subset \complex^N$ and a complex torus $K_\rho^\complex$ of dimension $N-n$ acting on $U_\FF$ so that we recover $X$ as the quotient:
\begin{equation}  
  X=U_\FF/K_\rho^\complex.
\end{equation}
The pair $(U_\FF,K_\rho^\complex)$ is defined as follows:
\begin{itemize}
\item For every $z=(z_1,\ldots, z_N)\in\complex^N$, we define $I_z:=\{j : z_j=0\} \subset \{0,\ldots,N\}$. Then we have that $$U_\FF := \{z : I_z \in \FF\}.$$ Clearly $U_\FF$ is the complement of a union of subspaces of $\complex^N$ of codimension at least 2. 
\item The complex torus $K_\rho^\complex$ acting on $\complex^N$ (and leaving $U_\FF$ invariant) is defined as the kernel of the the exact sequence: $$ 0 \to K_\rho^\complex \xrightarrow{i} \Torus^N \xrightarrow{\hat \rho} \Torus^n \to 0.$$ 
\end{itemize}
We will also use this sequence at the level of real Lie algebras: $$ 0 \to k_\rho \xrightarrow{i} \real^N \xrightarrow{\rho_\real} \real^n \to 0, $$ and its dual: $$0 \to \real^n \xrightarrow{\rho_\real^*} \real^N \xrightarrow{i_\rho^*} k^*_\rho \to 0.$$ Here we identify $\real^N$ with its own dual using  the canonical basis to induce the canonical inner product on $\real^n$.

It turns out that $K_\rho^\complex$ acts on $U_\FF$ with finite stabilizers when $P$ is a  Delzant polytope, and with no stabilizers at all when $P$ is an integral Delzant polytope, and therefore, $X$ is respectively an orbifold, or a manifold. For a proof of this fact we refer the reader to  proposition 2.2.3 page 157 of \cite{audin2004torus}.

\subsection{} Recall that a $2r$-dimensional \emph{symplectic manifold} $X$ is a smooth manifold together with a 2-form $\omega$ which is non-degenerate at each point, namely $\omega^{\wedge r}$ never vanishes. The non-degeneracy of $\omega$ tells us that it induces an isomorphism $T_x X \cong T^*_x X$ at every point. Given a smooth function $h : X \to \real$, we can obtain a vector field called its \emph{symplectic gradient} $\xi_h$ by looking at minus the image of $dh$ under the isomorphism $T^*X \to TX$, namely: $$\iota_{\xi_h} \omega = -dh,$$ giving a map: $$C^0(X) \to \mathrm{Vect} (X),$$ $$h \mapsto \xi_h.$$ Not every vector field is on the image of $\xi$. Given a vector field $\zeta$ of the form $\zeta = \xi_H$, we say that $H=H_\zeta$ is the \emph{Hamiltonian function} inducing the \emph{Hamiltonian vector field} $\zeta$. 

A group action of $G$ on $X$ induces by differentiation a Lie algebra homomorphism: $$\Psi:\mathrm{Lie}(G) \to \mathrm{Vect} (X).$$ The group action is called \emph{Hamiltonian} if every element in the image of this map is a Hamiltonian vector field. Whenever this happens for every element $a\in\mathrm{Lie}(G)$, we have a Hamiltonian function $H_{\Psi(a)}$ inducing $\Psi(a)$ so that $\xi_{H_{\Psi(a)}} = \Psi(a)$.

The \emph{moment map} of the Hamiltonian action of $G$ on $X$ is a map $\mu : X \to \mathrm{Lie}(G)^*$ defined by the following formula: $$\mu(x)(a) := H_{\Psi(a)}(x),$$ for every $x\in X$ and every $a \in \mathrm{Lie}(G)$.

For a regular value $\nu$ of the moment map $\mu$, the quotient $\mu^{-1}(\nu)/G$ makes sense because the Hamiltonian condition implies that $\mu^{-1}(\nu)$ is $G$-invariant with finite stabilizers, and that there is a unique symplectic form on the quotient lifting to the restriction of $\omega$ to $\mu^{-1}(\nu)$, making $\mu^{-1}(\nu)/G$ into a symplectic orbifold (cf. prop. 3.6.8 page 59 of \cite{audin2004torus}). In this case, the symplectic orbifold thus obtained $\mu^{-1}(\nu)/G$ is called the \emph{symplectic reduction} of the Hamiltonian action.

 \subsection{} The relation of the construction of a toric variety for a given Delzant polytope $P$ to the moment map can be easily described. In this section, a Delzant polytope $P\subset \real^n$ from which we compute $\FF$, $U_\FF$, $\rho$ and $K_\rho^\complex$ as before.
 
 The complex $N$-dimensional euclidean space $\complex^N$ is a symplectic manifold, for if we set its canonical coordinates to be $z_j = x_j + i y_j$, then the canonical symplectic form can be written as $\omega = \sum_j dx_j \wedge dy_j.$ The real torus $\mathbb{T}_\real^N$ has Lie algebra canonically isomorphic to $\real^N$ which we, in turn, identify with its own dual $(\real^N)^*$ as before. The action of the torus on euclidean space given by: $$(t_1,\ldots,t_N) \cdot (z_1,\ldots z_N) = (t_1 z_1, \ldots, t_N z_N)$$ turns out to be Hamiltonian, and a short calculation shows that its moment map $\tau \colon \complex^N \to \real^N$ is given by: $$\tau(z_1,\ldots,z_N):=(\pi|z_1|^2,\ldots,\pi|z_N|^2).$$ Similarly, the moment map $\mu_\rho \colon \complex^N \to k^*$ of the Hamiltonian action of $K_\rho^\real$ on $\complex^N$ is the composition: $$\mu_\rho \colon \complex^N \xrightarrow{\tau} \real^N \xrightarrow{i_\rho^*} k_\rho^*.$$ For a given Delzant polytope $P$ there is a $\nu_P \in k_\rho^*$ so that:
\begin{itemize}
\item The inverse image of the moment map defines a slice for $U_\FF$ with respect to $K_\rho^\complex$, namely: $$U_\FF = K^\complex_\rho \cdot \mu^{-1}(\nu_P).$$
\item The toric variety $X=X_P$ associated to $P$ is given by: $$ X_P = \mu^{-1}(\nu_P)/K^\real_\rho= U_\FF/K^\complex_\rho.$$ The vector $\nu_P$ contains all the information of the combinatorics $\FF$ of $P$.
\end{itemize}
To determine $\nu_P$ find real numbers $\lambda_1,\ldots,\lambda_N$ so that we can write the polytope as: $$P = \{ x\in \real^n \colon \langle x, \rho_j\rangle \geq \lambda_j,\ j=1,\ldots,N\},$$ then we have that $$\nu_P:=i_\rho^*(-\lambda_P),$$ where $\lambda_P \in \real^N$ is the vector $\lambda_P:=(\lambda_1,\ldots,\lambda_N).$






\subsection{}  We define a \emph{strictly convex polyhedral cone} generated by the $s$ vectors $v_1, \ldots ,v_s$ in $\real^n$ to be a set of the form: $$\sigma:=\{r_1 v_1+\cdots r_s v_s \colon r_j \geq 0\} \subset \real^n,$$ which does not contain a line. The cone $\sigma$ is called a \emph{rational} cone if the vectors $v_1, \ldots ,v_s$ in $\real^n$ can be chosen to have integer coordinates: $v_1, \ldots ,v_s$ in $\integer^n$ . A \emph{fan} $\Sigma$ in $\integer^n$ is a set of rational strictly convex polyhedral cones $\sigma$ so that each face of a cone in $\Sigma$ is also a cone in $\sigma$ and the intersection of two cones in $\Sigma$ is a face of each (see \cite{fulton1993introduction} section 1.4).

Presently we define the fan associated to a Delzant polytope $P$. The \emph{fan} $\Sigma_P$ is a collection of strictly convex polyhedral cones whose union is all of $\real^n$ and so that each cone is dual to a face of $P$. All the cones of fan $\Sigma_P$ are generated by a subset of the set of vectors $\{\rho_1,\ldots,\rho_N\}$.
 The cones in $\Sigma_P$ are exactly in one to one correspondance  with the sets in $\FF$, namely the cones: $$\sigma_I:=\{r_1 \rho_{i_1}+\cdots r_s \rho_{i_s} \colon r_j \geq 0\},$$ for the sets of indices $I=\{i_1\ldots,i_s\}\in\FF$.

We define now an equivalence relation on the set of all Delzant polytopes inside $\real^n$.  We say that two Delzant polytopes $P$,$P'$ are equivalent if and only if the corresponding fans $\Sigma_P$ and $\Sigma_{P'}$ associated to them are equal: $$[P]=[P'] \iff\Sigma_P = \Sigma_{P'},$$ Where we denote by $[P]$ the equivalence class of $P$. It is the same to have the equivalence class of $P$ as it is to have its fan $\Sigma_P$.

This equivalence relation has an alternative description. Two polytopes are equivalent if one is obtained from the other by a translation also we declare two polytopes $P$, $P'$ as equivalent whenever they are combinatorially equivalent and all the facets $F_j$ of $P$ are parallel to the corresponding facets $F_j'$ of $P'$. 

\subsection{} An important remark to make  is that, since the Delzant polytopes $P$ that we are considering satisfy a rationality condition, then the fans $\Sigma$ that we consider also satisfy this rationality condition: to wit, whenever $\sigma \in \Sigma$ is a top dimensional cone in $\Sigma$, it is generated by a $\rational$-basis  of $\rational^n$ composed by the integer vectors $\{ \rho_{i_1}, \ldots, \rho_{i_n} \}.$ The toric variety $X=X_\Sigma=X_P$ is smooth if and only if for every top dimensional cone $\sigma$ of $\Sigma$, the vectors $\{ \rho_{i_1}, \ldots, \rho_{i_n} \}$ actually form a $\integer$-basis  of $\integer^n$, this is proved in page 29 of \cite{fulton1993introduction}.

\subsection{} For example, consider $P$ to be a unit square on $\real^2$ with sides parallel to the coordinate axis. $P$ has $9$ faces (4 vertices of dimension 0, 4 edges of dimension 1 and one dense 2-dimensional face). We can compute $\Sigma_P$; it has 9 cones: in dimension 0, the origin; in dimension 1, the four rays at the coordinate axis; and the four 2-dimensional cones are the four quadrants. The fan can be considered as dual to the polytope $P$. 

\begin{figure}[htb]
\includegraphics[height=1.3in,width=3in,angle=0]{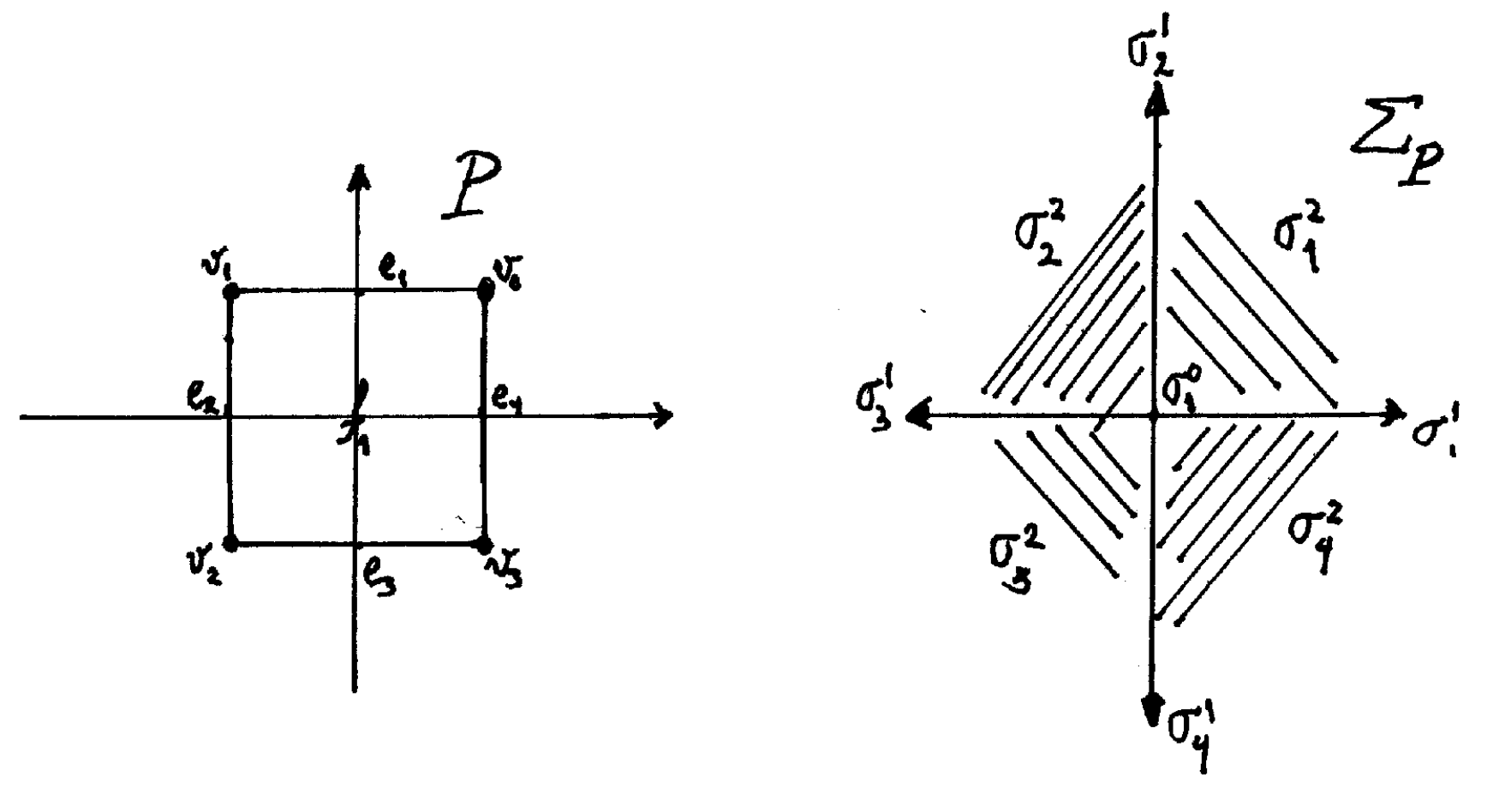}
\caption{A 2-dimensional polytope, and its fan.}
\label{firstfig}
\end{figure}

Observe that if $P'$ is a rectangle with sides parallel to the coordinate axis, we have that $[P]=[P']$ and $\Sigma_P = \Sigma_{P'}$, and as complex manifolds $X_P=X_{P'}$.

The toric variety $X_P$ is topologically the product of two 2-dimensional spheres; for a short computation shows that, by reordering the coordinates of $\complex^4$ ($N=4$ is the number of facets) in the manner $z=(z_1,z_3,z_2,z_4)$, we have that $U_\FF = (\complex^2  \setminus \{0\}) \times (\complex^2  \setminus \{0\})$ and $K^\complex_\rho = \complex^* \times \complex^*$, therefore:
$$ X_P = U_\FF/K^\complex_\rho = ((\complex^2  \setminus \{0\}) \times (\complex^2  \setminus \{0\}))/(\complex^* \times \complex^*) = \proj^1 \times \proj^1.$$

\subsection{} We leave to the reader to formulate the definition of  the concept of morphism of toric varieties $X\to X'$. Such maps should be algebraic, and also equivariant when restricted to the dense complex tori inside the varieties with respect to some homomorphism of the tori. In this way, we define the category $\Torics$ of toric varieties. It is not hard to define naturally a $\Fans$ category of (rational, Delzant) fans: its morphisms are the natural geometric combinatorial maps preserving integer lattices and sending cones to cones in an appropriate manner. Detailed definitions for the corresponding morphisms can be found, for example, in section 3, chapter 3 of \cite{CoxToric2011}.

The categories $\Torics$ and $\Fans$ are isomorphic (and, in particular, equivalent). The functor $\Torics \to \Fans$ goes by first assigning to $X$ the image of its moment map $P$ and to $P$ its associated cone $\Sigma_P$.

Let us describe the construction of the inverse functor $\Fans \to \Torics$. So to a fan $\Sigma$, we must associate a toric variety $X_\Sigma$ in a functorial way. It is enough to define the functor from  the full subcategory of cones $\Cones \to \Torics$ (but a warning here: the resulting toric variety will not be compact. Nevertheless, whenever the union of the cones of the fan $\Sigma$ cover all of $\real^n$, the resulting toric variety is compact.) Given a cone $\sigma \subset V:=\real^n$, we define the \emph{dual cone} by: $$\sigma^\vee:= \{ w \in V^*\cong \real^n \colon w(v) =  \langle v,w\rangle  \geq 0, \ \forall v\in \sigma \}.$$ The semigroup $S_\sigma$ associated to $\sigma$ is defined as $\sigma^\vee \cap \integer^n \subset \integer^n$, that $S_\sigma$ is finitely generated is the content of \emph{Gordon's Lemma}. We further associate a ring to $\sigma$: the semigroup ring $\complex[S_\sigma]$ of formal linear combinations of elements of $S_\sigma$ with coefficients in $\complex$, then the toric variety associated to $\sigma$ is: $$X_\sigma := \mathrm{Hom}_\mathrm{sg}(S_\sigma,\complex) = \mathrm{Specm}(\complex[S_\sigma]),$$ namely the points of $X_\sigma$ are the semigroup homomorphisms from $S$ to $\complex$ (cf. page 16 of \cite{fulton1993introduction}). The construction is clearly functorial, and thus can be extended to a functor $\Fans \to \Torics$. What is not obvious but still true is that this is the inverse functor to the moment map construction above.

One relevant remark is in order. The smallest non-trivial fan $\Sigma_0$ consists of only one 0-dimensional cone at the origin. Its associated toric variety is the complex torus, $X_{\Sigma_0} = \Torus^n$. Notice that we always have $\Sigma_0 \subset \Sigma_X$ and this immediately implies that  $\Torus^n \subset X$ as is required by the original definition.

It is possible to define a reasonable category $\Polytopes$ of equivalence classes of Delzant polytopes, and to prove that $\Polytopes^\op \cong \Fans \cong \Torics;$ the important point to remember is that $X$, $[P]$ and $\Sigma$ contain essentially the same amount of information.\footnote{More precisely, to have $P$ is the same as to have $X$ (as a complex manifold) plus a choice of a point in its Mori cone (K\"ahler cone).}

\subsection{} Given a fan $\Sigma$, let us write $X(\Sigma)$ to denote its associated toric variety.  Given fans $\Sigma'$ and $\Sigma$ we say that $\Sigma'$ is a refinement of $\Sigma$ whenever we have that each cone of $\Sigma$ is a union of cones in $\Sigma'$ and, in this situation, we have a canonical morphism of fans $\Sigma' \to \Sigma$ inducing the identity on lattices $\integer^n$. This, in turn, gives us a morphism $X(\Sigma') \to X(\Sigma)$ of the corresponding toric varieties which turns out to be birrational and proper.  We have that \emph{for any toric variety $X(\Sigma)$, there is a refinement $\Sigma'$ of $\Sigma$ so that $X(\Sigma') \to X(\Sigma)$ is a resolution of singularities}. This is proved in page 48 of \cite{fulton1993introduction}. 

For example, Delzant polygons $P\subset \real^2$ correspond to compact closed toric surfaces $X_P$ of real dimension 4. Their fans $\Sigma_P$ are completely determined by the primitive generators $v_1,\ldots,v_N \in \integer^2$ (in counterclockwise order) of its one dimensional cones (orthogonal to the $N$ faces of $P$). We think of the indices $i$ of the $v_i$ as elements of the cyclic group $i\in \integer/N\integer$. The variety $X_P$ is smooth only when $(v_i,v_{i+1})$ in a $\integer$-basis of $\integer^n$ for all $i$. As a matter of illustration consider a fan so that this condition fails only for $(v_1,v_2)$ where $v_2=(0,1)$ and $v_1=(m,-k)$ for $0<k<m$ with $k$ and $m$ relatively prime (cf. page45 of \cite{fulton1993introduction}). Write the \emph{Hirzeburch-Jung continued fraction} of $m/k$: $$\frac{m}{k}=a_1 -\cfrac{1}{a_2-\cfrac{1}{a_3-\cfrac{1}{\cdots-\cfrac{1}{a_r}}}}.$$
Then the fan $\Sigma'$, realizing the resolution of singularities for $\Sigma_P$, can be obtained by dividing $r$ times the cone generated by $v_1$ and $v_2$, where the slope of the $j$-th ray of said subdivision is completely determined by $a_j$. The process finishes because the continued fraction expansion is finite. For a far more detailed exposition of the relation of continued fractions and toric surfaces, we refer the reader to section 10.3 of \cite{CoxToric2011}.

\section{Non-commutative Geometry}

\subsection{} In this section, we very schematically describe the main structures appearing in the topological study of non-commutative spaces. We recommend the reader \cite{connesnoncommutative} for details.

,\subsection{} One motivation for non-commutative geometry is quantization; indeed, Heisenberg's uncertainty principle questions the notion of points and manifolds as good models for space-time. A second more mathematical formulation of the above is the \emph{commutative Gelfand-Naimark} theorem. It states that, for an arbitrary unital commutative $C^*$-algebra $A$, there is a compact Hausdorff topological space so that the algebra in question is isomorphic to the continuous real valued functions on $X$, namely $A\cong C^0(X)$ (in fact, $X$ can be canonically chosen to be the spectrum of $A$).\footnote{If we are to deal with non-compact, but still locally compact spaces, we need to consider non-unital algebras.} 

The Gelfand-Naimark theorem provides an equivalence of categories between $\commAlg^\op$, the opposite of the category of commutative $C^*$-algebras, and $\Top$, the category of Hausdorff locally compact topological spaces. Moreover, we have a pleasant dictionary between the algebraic properties of $A$ say being unital (resp. separable, having no projections) and $X$ correspondingly being compact (resp. metric, connected).

The dictionary can be extended to the usual concepts of rational algebraic topology: cohomology, $K$-theory and even rational homotopy theory. For example, when $X$ is a manifold the differential graded algebra that computes the cohomology  $H^*(X,\rational)$ is the \emph{Hochschild homology} ${HH}_\bullet(A)$ only depending on the algebraic structure of $A$ (we will return to its definition below).

\subsection{} We would like to fill the gap in the following diagram of categories: 
$$\begin{CD}
\Top     @>\cong>>  \commAlg^\op \\
@VV V        @VV V\\
?     @>\cong>>  \Alg^\op
\end{CD}$$

The most tautological way to resolve this problem is to define the category of non-commutative spaces $\NCTop$ as equal to the category $\Alg^\op$: a non-commutative space is just a non-commutative algebra. The non-commutative algebra $A$ is to be thought as the algebra of functions over the non-commutative space.

\subsection{} The rational algebraic topology of a non-commutative space can be written in terms of its non-commutatiuve algebra. Let $A$ be a unital, associative, possibly non-commutative algebra. We define the Hochschild
complex $C_\bullet (A,A)$ of $A$ as a negatively graded complex (for we prefer to have all differentials of degree $+1$):
$$ \stackrel{\partial}{\longrightarrow} A\otimes A\otimes A \otimes A \stackrel{\partial}{\longrightarrow}
A \otimes A \otimes A \stackrel{\partial}{\longrightarrow} A\otimes A \stackrel{\partial}{\longrightarrow} A,$$
where $A^{\otimes k}$ lives on degree $-k+1$.
The differential $\partial$ is given by
$$\partial(a_0 \otimes \cdots \otimes a_n) = a_0 a_1 \otimes a_2 \otimes \cdots \otimes a_n - a_0 \otimes a_1 a_2 \otimes \cdots \otimes a_n $$
$$+ \ldots + (-1)^{n-1} a_0 \otimes a_1 \otimes \cdots \otimes a_{n-1} a_n + (-1)^n a_n a_0 \otimes a_1 \otimes \cdots \otimes a_{n-1}.$$
This formula is more natural when we write the terms cyclically:
\begin{equation}
\begin{array}{ccccccc} &  &  & a_0 &  &  &  \\
&  & \otimes &  & \otimes &  &  \\  & a_n &  &  &  & a_1 &  \\  & \otimes
 &  &  &  &  \otimes & \\  & \vdots &  &  & &   \vdots  & \\ &  & \otimes & & \otimes &  & \\  & & & a_i & & &\end{array}
\end{equation}
for $a_0 \otimes \cdots \otimes a_n$. It is not hard to verify that $\partial^2=0$.

 The homology of the Hochschild complex has an homological algebra meaning:
 $$\mathrm{Ker\ }\partial/\mathrm{Im\ }\partial = \mathrm{Tor}_\bullet^{A\otimes_k A^{\mathrm{op}} - \mathrm{mod}}(A,A).$$

 An idea in non-commutative geometry (A. Connes) is that as $A$ replaces a commutative space the Hochschild homology of $A$ replaces, in
turn, the complex of differential forms.

\begin{theorem}[Hochschild-Konstant-Rosenberg, 1961, \cite{HKR}] Let $X$ be a smooth affine algebraic  variety, then if $A=\OO(X)$,
we have: $$HH_i(X):=H^{-i}(C_\bullet(A,A);\partial) \cong \Omega^i(X)$$ where $\Omega^i(X)$ is the space of $i$-forms on $X$.
\end{theorem}

To prove this, we consider the diagonal embedding $X\stackrel{\Delta}{\longrightarrow} X \times X$ and, by
remembering that the normal bundle of $\Delta$ is the tangent bundle of $X$, we have:
$$HH_\bullet(X) = \mathrm{Tor}_\bullet^{{\mathrm{Quasi-coherent}}(X\times X)} (\OO_\Delta,\OO_\Delta).$$ This, together with a
local calculation, gives the result.

 The Hochschild-Konstant-Rosenberg theorem motivates us to think of $HH_i(A)$ as a space of differential forms of degree $i$
on a non-commutative space.

Note that if $A$ is non-commutative, we have:
$$H^0(C_\bullet(A,A);\partial)=A/[A,A].$$ Also, for commutative
$A=\OO(X)$,  given an element $a_0 \otimes\cdots \otimes a_n$ in
$C_\bullet (A,A),$ the corresponding form is given by $\frac{1}{n!}
a_0 da_1\wedge \ldots \wedge da_n.$

 There is a reduced version of the complex $C^\redu_\bullet(A,A)$ with the same cohomology obtained by reducing modulo
constants, all but the first
factor: $$\longrightarrow A \otimes A/({\bf k}\cdot 1) \otimes A/({\bf k}\cdot 1)
\longrightarrow A \otimes A/({\bf k}\cdot 1) \longrightarrow A.$$

Alain Connes' main observation is that we can write a formula for an additional differential
$B$ on $C_\bullet(A,A)$ of degree $-1$, inducing a differential on
$HH_\bullet(A)$ that generalizes the de Rham differential:
$$B(a_0\otimes a_1 \otimes \cdots \otimes a_n) = \sum_\sigma (-1)^\sigma 1 \otimes a_{\sigma(0)}
\otimes \cdots \otimes a_{\sigma(n)}$$
where $\sigma\in \integer/(n+1)\integer$ runs over all cyclic permutations. It is easy to verify
that: $$B^2=0,\ \ \ B\partial +\partial B = 0, \ \ \  \partial^2=0,$$ which we depict pictorially as:
$$\xymatrix{
                   \cdots \ar@/_/[rr]_\partial && A \otimes A/1 \otimes A/1 \ar@/_/[ll]_B \ar@/_/[rr]_\partial &&
                   A \otimes A/1 \ar@/_/[ll]_B \ar@/_/[rr]_\partial && A \ar@/_/[ll]_B
                  }$$
 and  taking cohomology gives us a complex $(\Ker \partial / \Img
\partial ; B)$. A naive definition on the de Rham cohomology in
this context is the homology of this complex $\Ker B / \Img B.$

 We can do better \cite{kontsevich2008xi} by defining the negative cyclic complex $C^-_\bullet(A)$, which is formally a
projective limit (here $u$ is a formal variable, $\deg(u)=+2$):
$$C_\bullet ^-:= (C_\bullet^\redu(A,A)[[u]] ; \partial + u B) =
\lim_{\stackrel{\longleftarrow}{N}}(C_\bullet^\redu(A,A)[u]/u^N;\partial + u B).$$

 We define the periodic complex as an inductive limit:
$$C_\bullet ^\per:= (C_\bullet^\redu(A,A)((u)) ;
\partial + u B) = \lim_{\stackrel{\longrightarrow}{i}}(u^{-i} C_\bullet^\redu(A,A)[[u]];\partial + u B).$$
This turns out to be a ${\bf k}((u))$-module, this implies that
the multiplication by $u$ induces a sort of Bott periodicity. The
resulting cohomology groups called (even, odd) periodic cyclic
homology and are written (respectively): $$HP_\even(A), \ \ \ \ \
HP_\odd(A).$$ This is the desired replacement for de Rham
cohomology.

For example, when $A=C^\infty(X)$ is considered as a nuclear Fr\'echet algebra,
and if we interpret the symbol $\otimes$ as the topological tensor product then we have the canonical isomorphisms:
$$ HP_\even(A) \cong H^0(X,\complex) \oplus H^2(X,\complex) \oplus \cdots $$
$$ HP_\odd(A) \cong H^1(X,\complex) \oplus H^3(X,\complex) \oplus \cdots $$
\begin{theorem}[Feigin-Tsygan, \cite{Tsygan}] If $X$ is a affine algebraic variety ( \emph{possibly singular})
and $X_\topo$ its underlying topological space then:
$$HP_\even(A) \cong H^\even(X_\topo,\complex)$$ and $$HP_\odd(A) \cong H^\odd(X_\topo,\complex)$$
(these spaces
are finite-dimensional).
\end{theorem}

\subsection{} As we will see shortly, the most appropriate notion for a morphism in $\Alg$ between algebras $A$ and $B$ is that of a $A$-$B$-bimodule $M$. Given a $A$-$B$-bimodule $M$ and a $B$-$C$-bimodule $N$, we define the composition by using the tensor product: 
$$N\circ M := M \otimes_B N.$$ This gives a new notion of isomorphism of algebras. Two algebras $A$ and $B$ are called \emph{Morita equivalent} whenever there is an $A$-$B$-bimodule $M$, and a $B$-$A$-bimodule $N$ so that $M\otimes_B N \cong A$ (as $A$-$A$-bimodules), and $N\otimes_A M \cong B$ (as $B$-$B$-bimodules). It is easy to see that $A$ and $B$ are Morita equivalent if and only if their categories of modules, $A$-$\Mod$ and $B$-$\Mod$ are equivalent (and, therefore, the $K$-theories of $A$ and $B$ coincide).

A non-commutative algebra could be Morita equivalent to a commutative one. For example, it's very easy to check that the algebra of $n\times n$-matrices $\Mat_n(\complex)$ is Morita equivalent to the algebra $\complex$ for all $n\geq 1$.

Also if $A$ is Morita equivalent to $B$, then $HP_\bullet(A) \cong HP_\bullet(B)$.

\subsection{} Non-commutative spaces arise naturally from bad quotients $M/G$. If both $M$ and $G$ are compact and the action is free, the quotient naturally admits the structure of a smooth manifold and the map $M \to M/G$ is a submersion. But it is often the case that we have fixed points, or the Lie group $G$ acting on $M$ is not compact and a particular orbit may be dense on a open set of $M$. 

A standard procedure to solve this issue is to extend the category $\Man$ of smooth manifolds, embedding it fully faithfully into a bigger category, say: $$\Man \subset \Stacks.$$ A \emph{stack} is very much like an orbifold: it has local charts that look like $U/G$ for $G$ a compact Lie group acting smoothly on the manifold $U$. The charts are, as always, not canonical, and the group is not unique: it changes as we move on the stack, very much like in an orbifold. A stack is an orbifold if and only if all its local groups are finite \cite{MoerdijkSurveyOrbifolds, NearlyBook}. 

We write $[M/G]$ to indicate that the quotient is taken in the category of stacks. It is possible to choose a different manifold $M'$ and a different Lie group $G'$ so that $[M'/G']=[M/G]$, for example, take $M' := M \times \integer_2$, and $G':=G\times \integer_2$.

\subsection{} For example; consider the stack $\BB G := [\bullet/G]$, of $G$ acting trivially on a point $\bullet$. By the Yoneda lemma, the stack is completely determined by its functor of points $\Man^\op \to \Groupoids,$\footnote{Usually, for ordinary manifolds functors of points would be functors of the form $\Man^\op\to \Sets$; the whole point being that, for a stack, the functor must take values in $\Groupoids$.} from the category of manifold to the category of discrete groupoids, given by: $$M \mapsto C^\infty(M,\BB G):= \{ G \to P \to M \colon P\  \mbox{is a principal bundle over}\ M\}.$$ We shall interpret $C^\infty(M,\BB G)$ not as a set, but as a groupoid, the groupoid associated to the $G$-set of principal $G$-bundles over $M$.

Notice that the map $\bullet \to [\bullet/G] = \BB G$ becomes a submersion, and that the dimension of $\BB G$ is $- \dim G$. In fact the map, $\bullet \to \BB G$ is the universal principal $G$-bundle over $\BB G$; for given any $f_P=P \in C^\infty(M, \BB G),$ then $f_P^* \bullet = P$ tautologically:
$$\begin{CD}
P     @> >> \bullet \\
@VV V        @VV V\\
M     @>f_P >>  \BB G
\end{CD}$$

\subsection{} More generally, consider the stack $\XX:=[X/G]$. The functor of points $\Man \to \Groupoids$ assigns to a manifold $M$ the groupoid $\XX(M):=C^\infty(M,[X/G])$ whose objects are pairs $(P,\phi)$, with $P$ a principal $G$-bundle over $M$ and $\phi\colon M \to X$ a $G$-equivariant map. We say that $\XX$ classifies such pairs $(P,\phi)$. 

Let us consider the case $M=\bullet$. A point $(P,p_0) \in C^\infty(\bullet ,\XX) =\{ (P,\phi)\colon \phi \colon \bullet \to X, p_0:=\phi(\bullet) \}$ is a $G$-bundle $P=G$ over $X$, together with a choice of a marked point $p_0\in P$ (to have $p_0$ being the same as having $\phi$) . The automorphisms of such a point (the local stabilizer at $(P,p_0)$) is the reduced gauge group.

Notice that as $U$ ranges over all manifolds $C^\infty(U,\XX) = \XX(U)$ becomes a \emph{sheaf of groupoids} over $\Man$; often stacks are defined as sheafs of groupoids for this reason.

\subsection{} In practice, we can think of a general stack as an equivalence class $[\GG]$ of a Lie groupoid $\GG$ under Morita equivalence, the Lie groupoid taking the place of an atlas for the stack. Presently, we explain these concepts. We refer the reader to \cite{MoerdijkFoliation} for a more detailed exposition of these issues.

For example, the \emph{translation groupoid} $\GG = M \rtimes G$ is defined as the groupoid whose manifold of objects (homeomorphic to $M$) are the points $m\in M$, and whose manifold of arrows (homeomorphic to $M\times G$,)  are pairs $(m,g)$ for $m\in M$ and $g\in G$: the arrow $(m,g)$ send $m$ to $mg$. Clearly we have that $[M \rtimes G] = [M/G]$.

A \emph{groupoid} is a category $\GG$ where every arrow admits an inverse. We will denote by $\GG_0$ its set of objects, $\GG_1$ its set of arrows and:
$$\GG_1 \times_{s,t} \GG_1  \mathrel{\mathop{\longrightarrow}^{m}} \GG_1 \mathrel{\mathop{\longrightarrow}^{i}} \GG_1 \mathrel{\mathop{\rightleftarrows}^{s,t}_u} \GG_0,$$ where $m$ is the composition of arrows map, $i$ is the inverse, $s$ is the source, $t$ is the target and $u$ is the identity.
 
A \emph{Lie groupoid}, additionally, is a small groupoid where the set of all objects $\GG_0$ is a manifold as well as the set of all arrows $\GG_1$; furthermore, all the structure maps of the groupoid $\GG$ are smooth. 

The notion of Lie groupoid homomorphism (smooth functor of Lie groupoids) is natural and left to the reader. A Lie groupoid homomorphism $\phi_i : \HH_i \to \GG_i$, $i=0,1$, is said to be an \emph{essential equivalence} if 
\begin{itemize}
\item $\phi$ induces a surjective submersion $(y,g) \mapsto t(g)$ from $$\HH_0 \times_{\GG_0} \GG_1 =\{ (y,g) | \phi(y)=s(g) \}$$ onto $\HH_0$; and
\item $\phi$ induces a diffeomorphism $h \mapsto(s(h)\phi(h),t(h))$ from $\HH_1$ to the pullback $\HH_0 \times_{\GG_0} \GG_1 \times_{\GG_0} \HH_0.$
\end{itemize}
Two Lie groupoids $\GG'$ and $\GG$ are said to be \emph{Morita equivalent} if there are essential equivalences $ \GG \leftarrow \HH \rightarrow \GG'$ from a third Lie Groupoid $\HH$.

We will have that the stacks $[\GG']=[\GG]$ if and only if $\GG$ and $\GG'$ are \emph{Morita equivalent} (which more or less means that $\GG$ and $\GG'$ are equivalent as categories in a consistent $C^\infty$ fashion). 

\subsection{} For the study of non-commutative toric varieties we will be only interested in \'{e}tale groupoids that arise from foliated manifolds $(M,\FF)$. Let $q$ be the codimension of the foliation. The holonomy groupoid $\HH = \mathrm{Holo}(M,\FF)$ is the Lie groupoid so that $\HH_0 = M$, and two objects $x,y$ in $M$ are connected by an arrow if and only if they belong to the same leaf $L$. In such case, the arrows from $x$ to $y$ correspond to homotopy classes of paths lying on $L$ starting at $x$ and ending at $y$. 

A Lie groupoid $(\GG,m,i,s,t,u)$ is said to be \emph{\'{e}tale} if $m$, $i$, $s$, $t$ and $u$ are local diffeomorphisms. The Lie groupoid $\HH = \mathrm{Holo}(M,\FF)$  is always Morita equivalent to an \'{e}tale groupoid, for if we take an embedded $q$-dimensional transversal manifold $T$ to the foliation that hits each leaf at least once then the restricted groupoid $\HH|_T$ (where we only consider objects in $T \subset M$) is an \'{e}tale groupoid which is Morita equivalent to $\HH = \mathrm{Holo}(M,\FF)$ (cf. \cite{MoerdijkFoliation}).

\subsection{} Consider a stack $\XX=[\GG]$ where $\GG$ is \'{e}tale. Presently, we explain how to obtain a non-commutative space out of it. The non-commutative algebra $A_\GG$ associated to an \'{e}tale groupoid $\GG$ is the \emph{convolution algebra} of $\GG$; its elements are compactly supported smooth complex valued functions on the manifold of all arrows $\GG_1$ of $\GG$,  $f\colon \GG_1 \to \complex$. The product $f * g$ of two functions is the convolution product:
$$(f*g)(\alpha) = \sum_{\beta\circ\gamma=\alpha} f(\beta)g(\gamma).$$ Where the sum makes sense because it ranges over a space that is discrete because $\GG$ is \'{e}tale and finite because the functions are compactly supported. 

Suppose that $\GG' \simeq \GG$, then $\XX=[\GG]=[\GG']$ and the algebras $A_\GG$ and $A_{\GG'}$ are Morita equivalent, so the non-commutative space only depends on $\XX$ and not on $\GG$.

But it is still possible that two different stacks $\XX$ and $\XX'$ give the same non-commutative space. For instance, consider the 0-dimensional stacks $\XX_1=[\bullet/\integer]$ and $\XX_2=[\integer/\{1\}]=\integer$. The first one is connected, while the second has infinitely many components. The groupoids $\GG_1 = \bullet \rtimes \integer$ and $\GG_2 = \integer \rtimes \{1\}$ are not equivalent. But the Fourier transform $\FF\colon A_{\GG_1} \to A_{\GG_2}$ is an isomorphism and, therefore, a Morita equivalence. 

\subsection{} The homotopy type of a stack $\XX=[\GG]$ is defined to be the homotopy type of Segal's classifying space $B\GG$. In the case that $\GG=M\rtimes G$ then $B\GG \simeq M\times_G EG$ is the Borel construction for the group action. When the groupoid $\GG$ is \'{e}tale,  Crainic \cite{CrainicMoerdijk, crainic2004cech} has proved that localizing at units, we have a natural isomorphism: $$HP_*(A_\GG)_{(1)} \cong H_*(B\GG; \complex).$$
 
\subsection{Non-commutative Tori.} Take, for example, the Kronecker linear foliations $\FF_\theta$ on the torus $T^2= \real^2 /\integer^2$: Their universal cover are the foliations of $\real^2$ by all parallel lines of slope $\theta$. Whenever $\alpha$ is rational, all leaves are torus knots; but, when $\theta$ is irrational, all leaves are homemorphic to $\real$, and every leave is dense everywhere on $T^2$.

\begin{figure}[htb]
\includegraphics[height=1.3in,width=3in,angle=0]{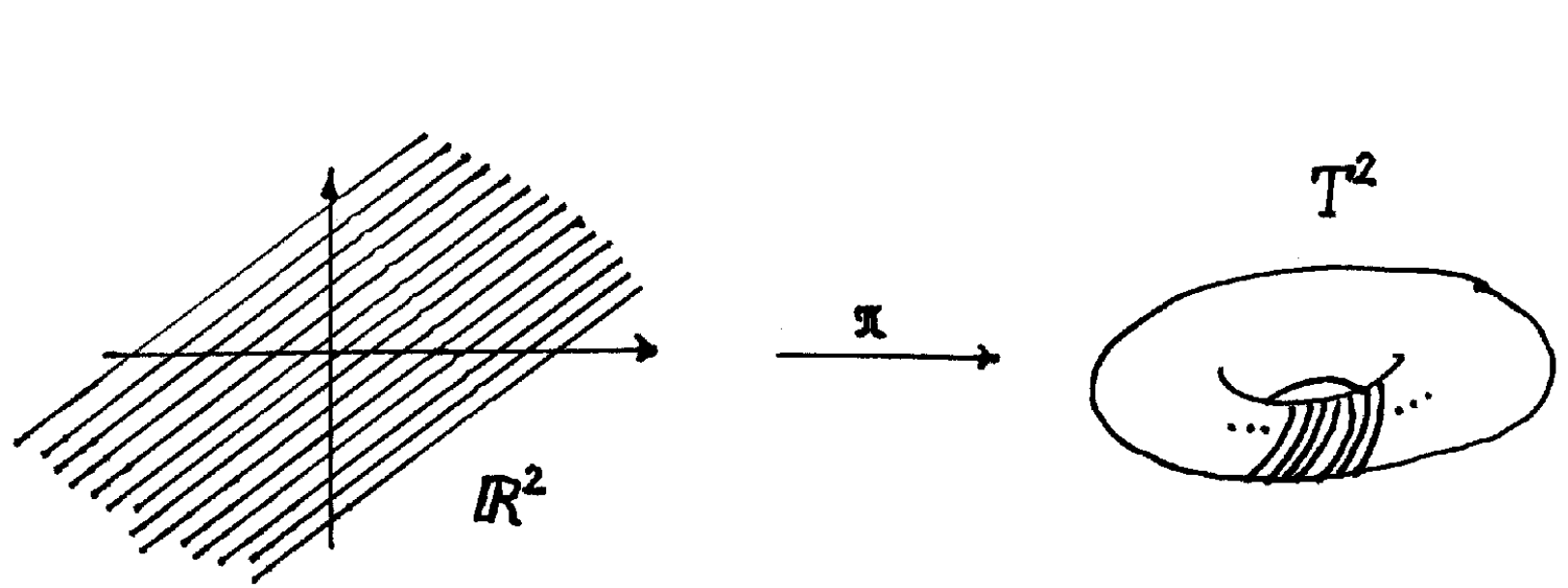}
\caption{The Kronecker foliation.}
\end{figure}

These foliations $\FF_\theta$ produce non-commutative algebras $A_\theta$ which are the convolution algebras of the holonomy groupoids associated to the foliations. Thought of as a non-commutative space, we write $T^2_\theta$ and we call it the \emph{non-commutative torus.} 

The algebra $A_\alpha \cong C^\infty(T^2_\theta) $ is after Fourier transform precisely the set of all
the expressions of the form: $$ \sum_{n,m\in\integer} a_{n,m} \hat{z}_1^n  \hat{z}_2^m, \ \ \ \ a_{n,m}\in\complex,$$ such
 that for all $k$, we have: $$a_{m,n} = O((1+|n|+|m|)^{-k}),$$ and $$\hat{z}_1 \hat{z}_2 = e^{i \theta} \hat{z}_2 \hat{z}_1.$$

 For $\theta \in 2\pi \integer$, we get the usual commutative torus.

A result of Marc Rieffel \cite{Rieffel} states that
$T^2_\theta$ is Morita equivalent to $T^2_{\theta'}$
if and only if: $$\theta' = \frac{a \theta +b}{c \theta + d},\ \ \ \
 \left(\begin{array}{cc}a & b \\c & d\end{array}\right)\in \mathrm{SL}_2(\integer).$$
And we have: $$ HP_\even(T^2_\theta)=H^0(T^2) \oplus H^2(T^2), $$ $$HP_\odd(T^2_\theta)=
H^1(T^2,\complex).$$

Finally, let us remark that, if we consider an elliptic curve: $$E=\complex/\left(\integer+ \tau \integer\right),\ \ \ \  \Im(\tau)>0.$$
Here: $$HP_\even(E) = H^0(T^2)\oplus H^2(T^2),$$ $$HP_\odd(E)=
H^1(T^2,\complex)$$

Non-commutative tori can be interpreted as 
limits of elliptic curves as $\tau\to \real$. 

\section{Variations on Diffeology}

\subsection{} In this section, we consider variations on a model for non-commutative spaces coming from non-singular foliations called \emph{diffeological spaces} \cite{iglesias2008diffeology}. We need to add additional structure to the classical diffeological spaces, and thus modify somewhat the classical definition.

\subsection{} Consider the category of topological spaces $\Top$. We will embed fully faithfully this category into the category $\Pro\Top$ of \emph{proto-topological spaces}; the category $\Pro\Top$ is a sort of Yoneda closure of $\Top$. 

We define a proto-topological space as a pair (\XX,X) of a functor\footnote{This is the functor of points of the proto-topological space. Notice that, unlike the case of stacks, functors of points go to $\Sets$ in this instance.} $\XX\colon \Top^\op \to \Sets$ (a pre-sheaf of sets on $\Top$), and a set $X:=\XX(\bullet)$ satisfying the following:
\begin{itemize}
\item For every topological space $Y$, we have: $\XX(Y) \subset \Hom_\Sets (Y,X).$
\item The set $\XX(Y)$ contains all the constant maps $Y\to X$ (therefore it is not empty).
\item The functor $\XX\colon\Top^\op \to\Sets$ is a sheaf; in particular, it sends the push-out diagram of inclusions:
$$\begin{CD}
U\cap V     @>{i_U}>> U \\
@VV{i_V}V        @VV{j_U}V\\
V     @>{j_V}>>  U \cup V
\end{CD}$$
to the pull-back diagram of restrictions:
$$\begin{CD}
\XX(U\cup V)     @>{\rho_U}>> \XX(U) \\
@VV{\rho_V}V        @VV{\rho_{U\cap V}}V\\
\XX(V)     @>{\rho_{U\cap V}}>>  \XX(U \cap V).
\end{CD}$$
\item The set $\XX(Y)$ admits a natural action by precomposition of $C^0(Z,Y)$, namely composition of maps induces a function: $$C^0(Z,Y) \times \XX(Y) \xrightarrow{\circ} \XX(Z).$$ 
\end{itemize}

The first two bullets above are part of the pre-sheaf structure but we spell them out explicitly to fix the notation.

We define $C^0(Y,\XX) := \XX(Y).$ Notice that, \emph{a priori}, $X$ does not carry a topology. 

The Yoneda lemma implies that $\Top \subset \Pro\Top$, by sending the topological space $X$ to the functor $\XX(Y):=C^0(Y,X).$

\subsection{} We can define the homotopy type of a proto-topological space $\XX$ by first associating to it  a simplicial set: $$X_n := \XX(\Delta^n),$$ where $\Delta^n$ is the standard $n$-simplex. Clearly $X_\bullet$ inherits a natural structure of a simplicial set, and its geometric realization $||X_\bullet||$ will be called the homotopy type $|\XX|$ of $\XX$.

Say, for instance, that we started by taking a \emph{bona fide} topological space $X$ (so that $\XX(Y):=C^0(Y,X),$) then it is a classical result from basic topology that there is a weak homotopy equivalence $X\simeq ||X_\bullet||$.

In any case, we define the cohomology $H^*(\XX)$ of $\XX$ to be $H^*(||X_\bullet||).$ The de Rham theorem is not true in general.

\subsection{} A \emph{diffeological space}  \cite{souriau1980groupes} is a pair $(\XX,X)$ of a functor $\XX\colon \Man^\op \to \Sets$ and a set $X:=\XX(\bullet)$ so that all the axioms of a proto-topological space are satisfied if we replace $\Top$ by $\Man$, and we change $C^0(Z,Y)$ by $C^\infty(Z,Y)$, the category of manifolds embeds fully faithfully into the category of diffeological spaces: 
$$\Man \subset \Diffeologies.$$

We can define differential forms on a diffeological space as follows. We say that a family of $k$-differential forms $\omega = \{ \omega_{f,M} \}$ indexed by pairs $(f,M)$ where $M \in \Man$ ranges throughout all manifolds, and $f\in \XX(M)$ is a $k$-differential form on $\XX$ if it satisfies the differential equation: 
$$ \omega_{f \circ g, N} = g^* \omega_{f,M}, $$ whenever $g\in C^\infty(N,M)$ and $f\in \XX(M)$.

Since Cartan's exterior differential $d$ commutes with pullbacks, this defines a dga $(\Omega^*(\XX),d)$ canonically associated to the diffeological space $\XX$. The cohomology algebra of this dga is called the \emph{de Rham cohomology} $H^*_\dR (\XX)$ of $\XX$.

\subsection{} Effective orbifolds (orbifolds where all the local groups act effectively) can be modelled by diffeological spaces.

Notice that it is enough to have the values $\XX(U)$ of $\XX$ at all open sets $U$ of euclidean spaces, for manifolds can be obtained gluing those. A map $p \in \XX(U)$ is called a \emph{plot} of the diffeology $\XX$.

It is an interesting fact that not every diffeology $\XX$ is of the form $\XX(Y)=C^\infty(Y,X)$ for some manifold $X$.

 It is the same to say that $ p \in \XX(U)$ than to say that $p: U\to X$ is a plot on $X$. Observe that the full functor $\XX$ can be recovered by having the set $X$ and  the prescription to decide when a map of sets $U \to X$ is a plot.

Iglesias, Karshon and Zadka have proposed a definition for an effective orbifold in terms of diffeologies \cite{iglesias2010orbifolds}. The following definitions are theirs:

Let $\XX$ be a diffeological space, let $\sim$ be an equivalence relation
on $X$, and let $\pi \colon X \to Y:=X/\!\sim$ be the quotient map.
The \emph{quotient diffeology} on $Y$
is the diffeology in which ${p}\colon{U}\to{Y}$ is a plot
if and only if each point in $U$ has a neighborhood $V \subset U$
and a plot ${\tilde{p}}:{V}\to{X}$
such that $p|_{V} = \pi \circ \Tilde{p}$.

A diffeological space $\XX$ is \emph{locally diffeomorphic}
to a diffeological space $\YY$ at a point $x \in X$ if and only if
there exists a subset $A$ of $X$, containing $x$, and there exists
a one-to-one function $f \colon A \to Y$ such that:
\begin{enumerate}
\item
for any plot $p \colon U \to X$,
the composition $f \circ p$ is a plot of $Y$;
\item
for any plot $q \colon V \to Y$, the composition
$f^{-1} \circ q$ is a plot of $X$.
\end{enumerate}

An $n$ dimensional manifold can be interpreted as a diffeological space
which is locally diffeomorphic to $\real^n$ at each point.

A \emph{diffeological orbifold} is a diffeological space
which is locally diffeomorphic at each point to a quotient
$\real^n/\Gamma$, for some $n$, where $\Gamma$ is a finite group
acting linearly on $\real^n$.

As expected, diffeological orbifolds with differentiable maps
form a subcategory of the category of diffeological spaces. Moreover, for effective orbifolds, one can recover the orbifold from the diffeological orbifold \cite{Karshon2008}.

Unfortunately, this does not work so well for orbifolds that are non-effective as, for example, $\XX:=[\bullet/G].$ The category of diffeological orbifolds contains the category of manifolds, but the category of finite groups does not fit nicely on this approach. To deal with non-effective orbifolds, we must think of functors: $$\XX\colon \Man^\op \to \Groupoids.$$ rather than of functors: $\XX \colon \Man^\op \to \Sets.$

\subsection{} We define a \emph{proto-manifold} as a pair $(\bar\XX,\XX)$ where $\bar\XX$ is a diffeological space and $\XX$ is a proto-topological space satisfying: $$ \bar\XX(M) = C^\infty(M,\bar\XX) \subset C^0(M,\XX) = \XX(M) $$ for every manifold $M$. Proto-manifolds have homotopy types defined by $|\bar\XX|:=|\XX| = || X_\bullet ||.$

A complex diffeological space is a functor $\tilde\XX \colon \Man_\complex^\op \to \Sets$ from the category of complex manifolds $\Man_\complex$ to the category of sets $\Sets$ satisfying the same axioms as those of a proto-topological space replacing $\Top$ by $\Man_\complex$, and changing $C^0(Z,Y)$ by $\Hol(Z,Y)$.

A proto-complex manifold is a triple $(\tilde\XX, \bar\XX,\XX)$ so that: $$\Hol(M,\tilde\XX) \subset C^\infty(M,\bar\XX) \subset C^0(M,\XX)$$ for every complex manifold $M$. By abuse of notation, we will simply write $\XX$ for $\tilde\XX$, $\bar\XX$ and $\XX$.

We could also define the concept of K\"ahler diffeology, but as it gets slightly technical, we will refer the reader to \cite{NCTorics} for details.  

It is also possible to define the moment map and its image for  diffeological spaces \cite{iglesias2008diffeology}.

\subsection{} Given a diffeological space (resp. proto-topological space, proto-complex manifold) $\XX$ and an equivalence relation $\sim$ on $X$, we can endow $X/\sim$ with the structure of a diffeological space (resp. proto-topological space, proto-complex manifold) canonically as follows. Let $\pi\colon X \to X/\sim$ be the canonical projection. For an arbitrary manifold (resp. topological space, complex manifold), we define the functor $\XX/\sim$ by taking the smallest diffeology\footnote{This is also called the final diffeology associated to the map $\pi$ (cf. \cite{HectorLie}) and, as a referee has kindly pointed out, corresponds to a sheafification of the pre- sheaf $ \DD(M) =  \{ \pi\circ f \colon f\in \XX(M)\} $.} so that:  $$ \{ \pi\circ f \colon f\in \XX(U)\} \cup \{\mathrm{constants}\} \subseteq (\XX/\sim) (U)$$ for every open domain $U$ of an euclidean space. Thus we have arbitrary quotients in the categories $\Pro\Top$, $\Diffeologies$ and $\Pro\Man_\complex$.

A very important sub-example is that of the space of leaves of a foliation. Let $\FF$ be a (holomorphic) foliation on a (complex) manifold $X$, and let $\XX$ be the (proto-complex manifold) diffeology associated to $X$ by $\XX(M):=C^\infty(M,X)$ (resp. $\XX(M) = \Hol(M,X)$). Then, by taking $\sim$ to be the equivalence relation on $X$ given by $x\sim y \iff x\in L$ and $y\in L$ for the same leave $L$ of  $\FF$, we can consider the quotient set $X/\FF$ called the space of leaves of $\FF$. It acquires canonically the structure of a (proto-complex manifold) diffeological space $\XX/\FF$.

To relate this to the previous discussions, consider the holonomy groupoid $\GG$ of $\FF$ and let's point out that it is quite easy to verify that, as (proto-complex manifolds), diffeological spaces $\GG_1/\GG_0$ and $\XX/\FF$ are isomorphic.

It is known \cite{hector2011rham, crainic2004cech} that the de Rham cohomology $H^*_\dR (\XX/\FF)$ is isomorphic to the base-like cohomology $H^*_b(X,\FF)$ and, in turn, isomorphic to the localization at units of periodic cyclic cohomology $HP^*(A_\GG)_{(1)}$. This is an important relation of the theory of diffeologies with non-commutative geometry.

\subsection{} In particular, think again of the non-commutative torus $\TT^2_\theta$, but this time as a diffeological space: As the  space of leaves $T^2/\FF_\theta$ of the Kronecker foliation with the quotient diffeological structure. It is a fun exercise to show that $\TT^2_\theta \cong \TT^2_{\theta'}$ as diffeological spaces if and only if: 
$$\theta' = \frac{a \theta +b}{c \theta + d},\ \ \ \
 \left(\begin{array}{cc}a & b \\c & d\end{array}\right)\in \mathrm{SL}_2(\integer),$$ and by the theorem of Rieffel, this occurs if and only if $A_\theta$ is Morita equivalent to $A_{\theta'}$. This is the second important relation that we encounter between diffeologies and non-commutative geometry.

\section{LVM-theory}

\subsection{} The theory of LVM-manifolds and their associated structures started in 1984 at a seminar in Mexico City where the fourth author of this paper proposed the study of the topology  of the intersection of real quadrics $M$ of the form: $$\sum_{i=1}^n \lambda_i |z_i|^2 =0,$$ for $n$ fixed complex numbers $\lambda_1, \ldots \lambda_n$ to Lopez de Medrano. We refer the reader to \cite{de1997new, meersseman2004holomorphic} for full details. In this section, we will describe the aspects of the theory that we need most for the definition of non-commutative toric varieties.

\subsection{} The intersection of real quadrics $M$ defined above appears naturally when considering the system of linear, complex, differential equations in $\complex^n$: $$\dot z_i = \lambda z_i.$$
The solutions of this system define an action of $G=\complex$ on $X=\complex^n$, whose orbits foliate $\complex^n \setminus \{0\}$. A leaf $L$ is called a \emph{Poincar\'e leaf} if the origin belongs to its closure $\bar L$, and otherwise it is called a Siegel leaf. We will concentrate on the Siegel leaves. 
We will need to assume throughout that $\Lambda:=(\lambda_1,\ldots,\lambda_n) \in \complex^n$ satisfies the \emph{weak hyperbolicity condition}; namely we will assume that $0\in\complex$ is not contained in any of the segments $[\lambda_i,\lambda_j]$.
\begin{figure}[htb]
\includegraphics[height=1in,width=1in,angle=0]{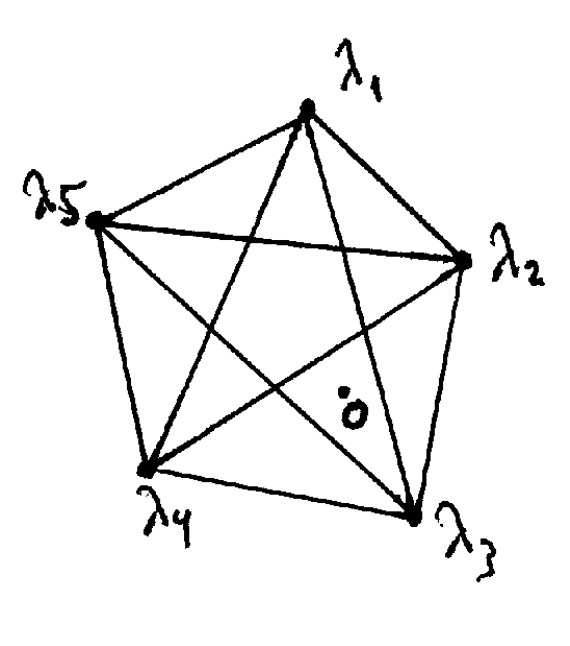}
\caption{The weak hyperbolicity condition.}
\end{figure}

We will further assume that $0$ is in the convex hull of the $\lambda_i$ as in the picture. 

We let $\SSS$ be the union of the Siegel leaves: With our assumptions, $\SSS$ is open and dense in $\complex^n$. The complement of $\SSS$ is the union of coordinate subspaces of $\complex^n$ for index sets so that the corresponding sub-configuration of $\Lambda$ does not contain $0$ in its convex hull. So $\Lambda$ already contains the combinatorics of the situation.

The (complex) time translation action of $G=\complex$ on $X=\complex^n$ clearly leaves $\SSS$ invariant, and commutes with the scalar multiplication action of $\complex^*$ on $\SSS$. The manifold $$N=N(\Lambda) := \SSS / \complex \times \complex^*,$$ is called an \emph{LV-manifold}. The manifold $N$ is related to the manifold $M$ as follows: $\SSS/\complex$ can be identified with a transversal manifold to the foliation, which can be easily computed by noticing that every leave $L$ inside $\SSS$ has a unique point closest to the origin, and hence: $$M:=\{ z\in \complex^n \setminus \{0\} \colon \sum \lambda_i |z_i|^2 =0 \} \cong \SSS/\complex,$$ is the transversal to the foliation restricted to $\SSS$. From this, we conclude that its projectivization $M/\complex^*$ satisfies: $$ N \cong M/\complex^* \subset \complex\proj^{n-1},$$ therefore, $N$ is the intersection of real quadrics inside $\complex\proj^{n-1}$.

For $n>3$ it can be shown that the LV-manifolds $N$ are smooth, compact, complex manifolds that \emph{do not admit a symplectic structure}. Topologically they are connected sums of products of spheres \cite{de1997new}.

\subsection{} More generally, we let $m$ and $n$ be two integers with $n>2m$, and we let $\Lambda:=(\Lambda_1,\ldots,\Lambda_n)$ be a set of $n$ vectors on $\complex^m$. We ask for this data to satisfy the \emph{Siegel condition}, that is to say that $0$ is in the real convex hull of $\Lambda$; and the \emph{weak hyperbolicity condition} that says that, for every subset of $\Lambda$ with $2m$ elements, $0$ is not in the convex hull of the subset $(\Lambda_{i_1},\ldots,\Lambda_{i_{2m}})$.  We then call $\Lambda$ admissible. With this data, we produce a dynamical system that amounts to a $G=\complex^m$ action on $X=\complex^n$, explicitly:
$$ (T,z) \in \complex^m \times \complex^n \mapsto (z_1 e^{\langle \Lambda_1, T \rangle},\ldots, z_n e^{\langle \Lambda_n, T \rangle})\in\complex^m.$$
Again as before, the projectivization of the space of leaves of the foliation restricted to the Siegel leaves is a complex manifold:
$$ N := \SSS/\complex^m \times \complex^* \cong \left\{ {[z]\in\complex\proj^{n-1} \colon \sum_{i=1}^n \Lambda_i |z_i|^2 =0} \right\}.$$ We call these manifolds LVM-manifolds.

Every diagonal holomorphic vector field $\xi = \sum_{i=1}^n \alpha_i z_i \frac{\partial}{\partial z_i}$ on $\complex^n$ descends onto a holomorphic vector field on $N$. In particular, if we write $\Lambda_i = (\lambda_i^1,\ldots\lambda_i^m),$ for $1\leq i \leq n$, and using the real parts $\Re(\lambda_j^i)$, we can define $m$ commuting holomorphic vector fields on $\SSS$:
$$\eta_i := \Re(\lambda_j^i) z_j \frac{\partial}{\partial z_j}, $$ that together define a locally free holomorphic action of $\complex^m$ on $N$ and, with this, a foliation $\FF$ on $N$ with leaves of dimension $m$. We call the pair $(N,\FF)$ an LVM-pair.

\subsection{} We say that the matrix $\Lambda$ satisfies the \emph{rationality condition} (K) if the space of solutions of the system (S) of equations:
$$\sum_{i=1}^n s_i \Lambda_i =0,$$
$$\sum_{i=1}^n  s_i =0,$$
admits a basis with integer coordinates.

We define the $n$-tuple of vectors $(v_1,\ldots,v_n)$ in $\real^{n-2m-1}$ by declaring that $(x_1,\ldots,x_n)$ is a solution of the system (S) if and only if $x_1=\langle v_1,u \rangle,\ldots, x_n=\langle v_n,u\rangle$ for some $u \in \real^{n-2m-1}$. The $n$-tuple of vectors $(v_1,\ldots,v_n)$ is called the \emph{Gale transform} of $\Lambda$.

The third and fourth authors of this paper proved:

\begin{theorem} 
Let $\Lambda$ be an admissible configuration satisfying condition (K). Then:
\begin{itemize}
\item The leaves of the foliation $\FF$ of $N$ are compact complex tori of dimension $m$.
\item The quotient space $N/\FF$ is a projective toric variety $X_P$ of dimension $n-2m-1$.
\item The LVM-manifold (which is a moment angle complex) is the 2-fold homotopy cover of $X_P$, namely $N\simeq X_P\langle 2 \rangle$.
\item The image $P$ of the moment map for $X_P$ is a convex polytope that can be written as the set of vectors $u\in\real^{n-2m-1}$, so that 
$$ \langle v_1,u \rangle \geq -\epsilon_1 ,\ldots,  \langle v_n,u\rangle \geq -\epsilon_n,$$ for certain real values $\epsilon_1\ldots \epsilon_n$.
\end{itemize}
Moreover, whenever $\Lambda$ does not satisfy (K) then $\FF$ has non-compact leaves dense on an open set of $N$.
\end{theorem}

Notice that the dichotomy that we had for the non-commutative torus between satisfying a rationality condition if and only if the leaves of the Kronecker foliation were compact and non-rationality implying non-compactness for the leaves has a perfect analogue here. 

The topology of $N\simeq X_P\langle 2 \rangle$ is that of a moment angle complex \cite{BP2002, BBCG2010}.

Finally, let us observe that, roughly speaking, $[P]$ and $\Lambda$ are in a duality correspondance given by the Gale transform; therefore, every $X_P$ is of the form $N(\Lambda)/\FF$ for $\Lambda$ the Gale transform of $P$.

\section{Non-commutative Toric Varieties}

\subsection{} In this section, we define non-commutative toric varieties and describe some of their most basic properties. We refer the reader to \cite{NCTorics} for details.

\subsection{} Let $\Lambda$ be an admissible configuration that (possibly) does not satisfy condition (K), we make the following definition:
\begin{definition}
The \emph{non-commutative toric variety associated to} $\Lambda$ is the proto-complex manifold: $$\XX(\Lambda):=[N(\Lambda)/\FF] = \NN(\Lambda)/\complex^m.$$ 
\end{definition}
If you forget the complex structures, then $\XX(\Lambda)$ can be seen as a diffeological space.

\subsection{} It is possible to show that non-commutative toric varieties are K\"ahler and, therefore, symplectic. They have a (proto-topological) moment map: 
$$\mu\colon \XX(\Lambda) \to P \subset \real^{n-2m-1}$$
and $P$ is the Gale dual to $\Lambda$, and so $\XX$ can be completely recovered from $P$, we will write $\XX_P := \XX(\Lambda)$.

All compact convex polytopes $P$ that appear from the previous construction are characterized by the following property: \emph{at each vertex $v$ of $P$, there are exactly $n$ edges, and there are vectors $(v_1,\ldots,v_n)$ along each of the $n$ edges meeting at $v$, forming a $\real$-basis of the space $\real^n$}. These polytopes will be called \emph{irrational Delzant polytopes}\footnote{In combinatorial geometry, a \emph{simple polytope} is defined as a $d$-dimensional polytope each of whose vertices are adjacent to exactly $d$ edges. In our terminology, a convex simple polytope is either an irrational or a rational Delzant polytope.}. The polytopes satisfying the rationality condition (at each vertex $v$ of $P$, there are exactly $n$ edges, and there exist the shortest integer vectors $(v_1,\ldots,v_n)$ along each of the $n$ edges meeting at $v$, forming a $\rational$-basis of the subgroup $\rational^n \subset \real^n$,) are called  \emph{(rational) Delzant polytopes} and the polytope $P$ is rational Delzant if and only if $\XX$ is an orbifold diffeology. The diffeology $\XX$ is a manifold diffeology (and not only an orbifold diffeology) if and only if, additionally at each $v$ the vectors $(v_1,\ldots,v_n)$, form a $\integer$-basis of the lattice $\integer^n \subset \rational^n$ (\emph{integral Delzant polytopes}).

Here we are lying a little bit. Consider the smooth integral case. In the case in which $P$ is integral Delzant, $\XX_P$ has the homotopy type of $N$ and not the homotopy type of $X_P=N/\FF=N/\complex^m$. What happens is that we must consider not $X_P$ but rather a pair $(X_P, \GG_P)$ where $\GG_P$ is a canonical gerbe\footnote{This is like saying that we must interpret $X_P$ as an orbifold with a stabilizer $\integer^m$ at every point.} with abelian band $\integer^m$ over $X_P$. The reason for this is that while, as topological spaces $X_P=N/\Torus^m = N/\complex^m$, nevertheless as diffeological spaces (as stacks, as non-commutative spaces): $$ [N/\Torus^m] \neq [N/\complex^m]. $$ What we really have is  $[N/\complex^m] = (X_P, \GG_P)$, for the $\complex^m$-action on $N$ has a  stabilizer $\integer^m$ at every point. Since, from the point of view of the classifying space $B\integer^m \simeq \Torus^m$, it is easy to see, using a spectral sequence, that the homotopy type of $\XX_P$ is that of the 2-fold homotopy cover $N$ of $X_P$ (a moment angle complex).

For example, if we define $\complex\proj^1$ as $[S^3/S^1]$ it then has the homotopy type of $S^2$ (by the Hopf fibration); but if we define $\widetilde{\complex\proj^1}:=[S^3/\real]$, then it is topologically an $S^2$ together with a gerbe with band $\integer$ and Dixmier-Douady characteristic class equal to the fundamental class of $S^2$; therefore, $\widetilde{\complex\proj^1}$ has homotopy type equal to $S^3$.

In this theory, when we say a non-commutative $\complex\proj^1$ we mean a deformation of $\widetilde{\complex\proj^1}$ and not a deformation of $\complex\proj^1$.

Another example is provided by the non-commutative torus: In the case $\theta$ rational, if we define $T^2_\theta$ as $T^2/S^1_\theta$ then $T^2_\theta$ would have the homotopy type of a circle, and we do not want that. But if we define $\TT^2_\theta$ as $[T^2/\real_\theta]$, every point acquires a stabilizer $\integer$ and the result is the circle $T^2/S^1_\theta$ together with a $\integer$-gerbe, hence the homotopy type of $\TT^2_\theta$ is $T^2$ as desired.

From now on, whenever we say toric variety (orbifold or smooth), we will mean $\widetilde{X_P}$ rather than $X_P$, but we will drop the tilde  and just write $X_P$. So the homotopy type of $X_P$ is the moment angle complex $N$. 
\subsection{} Deformation theory works for non-commutative toric varieties and allows us to prove:
\begin{theorem} Fixing the combinatorial structure $\FF$ of $P$, there exists moduli space $\MM_\FF$ of non-commutative toric varieties of the form $X_P$, and it is a complex orbifold. The rational points of $\MM_\FF$ correspond to ordinary classical toric varieties.
\end{theorem}
The point being, that the Gale transform gives us complex parameters $\Lambda$ from real parameters $P$.

\subsection{} There is a category of (irrational) fans and an equivalence of categories with non-commutative toric varieties. In other words, the theory of fans extends to the new case defining functorially $\XX_\Sigma$ for an irrational fan $\Sigma$.

\subsection{} Let us consider the example of non-commutative $\complex\proj^1$. For this, we need to take $n=4$, $m=1$ and the admissible configuration $\lambda_1=1=\lambda_3$, $\lambda_2=i$ and $\lambda_4 = -1-i$.

In this case, $N$ is a Hopf surface (homeomorphic to $S^3\times S^1$) fibering over $\complex\proj^1$ with fibers elliptic curves $E_\tau$ and, in this case, the leaf space of the foliation is $\complex\proj^1 = N/E_\tau$, $\widetilde{\complex\proj^1} = N/\complex$.

Let us vary $\lambda_3$. As $\lambda_3 = p = 3l+1$, $l>0$, we get that $\XX_p=[N_p/\FF]$ becomes a teardrop orbifold with a singularity of order $\frac{2p+1}{3}$. If $p$ is not equal to $1$ modulo 3, then $\XX_p$ is a football orbifold with two singularities, one of order 3 and a second one of order $2p+1$. 

When $\lambda_3$ is an irrational real number, $\XX_{\lambda_3}$ admits a moment map: $$\mu \colon  \XX_{\lambda_3} \to [0,1],$$ so that, if we write the composition $\nu:N\to \XX_{\lambda_3} \to [0,1],$ then $\nu^{-1}(0)$ and $\nu^{-1}(1)$ are leaves of the foliation homeomorphic to elliptic curves but, at any other, point $0<t<1$ the inverse image $\nu^{-1}(t)$ is a complex torus $\Torus^2$ foliated by an irrational slope $\theta_t$ by leaves homeomorphic to $\complex$. Therefore, $\mu^{-1}(0)$ and $\mu^{-1}(1)$ are points (with a $\integer \times \integer$ stabilizer each) while $\mu^{-1}(t)$ is a non-commutative complex torus.

\subsection{} We will have in general that, at a generic point $p\in P$, the inverse image $\mu^{-1}(p)$ of the moment map: $$\mu\colon \XX(\Lambda) \to P \subset \real^{n-2m-1},$$ will be a non-commutative torus (sometimes, of course, it will be an ordinary commutative torus). It will always be a commutative torus in the case in which $\Lambda$ satisfies condition (K).

\subsection{} Let us consider a second example: Take $n=5$, $m=1$ and the admissible configuration $\lambda_1=\lambda_4=1$, $\lambda_2=\lambda_3=i$ and $\lambda_5=2-2i$.

In this case, $N$ is a 3-fold (6 real dimensions) homeomorphic to $S^3 \times S^3$, for it is the intersection of two real quadrics inside $\complex\proj^4$ given by:
$$ |z_1|^2 + |z_4|^2 - 2|z_5|^2 = 0,$$
$$ |z_2|^2 + |z_3|^2 - 2|z_5|^2 = 0.$$
These equations can be rewritten as:
$$ 2|z_5|^2 = |z_1|^2 + |z_4|^2 = |z_2|^2 + |z_3|^2. $$
The conditions that define the Siegel domain are $z_5 \neq 0$, $(z_2,z_3) \neq (0,0)$ and $(z_4,z_1)\neq (0,0)$. So projectivizing (setting $z_5=1$), we have:
$$ 2 = |z_1|^2 + |z_4|^2 = |z_2|^2 + |z_3|^2, $$ proving that $N=S^3\times S^3$.

The composition $\nu \colon N \to X(\Lambda) \xrightarrow{\mu} \real^5$ is given by: $$\nu(z) = (|z_1|^2,\ldots, |z_5|^2).$$ From this, it is easy to calculate $\nu^{-1}(1,1,1,1,1) = S^1 \times S^1 \times S^1 \times S^1 \subset N.$

The manifold $X$ is $\complex\proj^1 \times \complex\proj^1$ and $\mu^{-1}(1,1) = S^1 \times S^1 \subset \complex\proj^1 \times \complex\proj^1.$ The quotient map $N\to X = \complex\proj^1 \times \complex\proj^1$ can be given in coordinates as: $$[z_2:z_3:z_4:z_1] \mapsto ([z_2:z_3],[z_4:z_1]).$$ which, at the level of fibers of the moment map, looks like two copies of the map $S^1\times S^1 \to S^1$ given by $(z_1,z_2) \mapsto z_1/z_2$, which in exponential form reads: $(e^{2\pi i \alpha}, e^{2\pi i \beta}) \mapsto e^{2\pi i (\alpha-\beta)}$. This, in turn, looks like $(\alpha, \beta) \mapsto \alpha-\beta$ as a map $\real \times \real \to \real$. The fibers of this map define the Kronecker foliation of slope $\theta=1$.

As we vary $\lambda_1$ and we make it $\lambda_1 = 1+ \epsilon$, $N$ becomes: 
$$ 2 = (1+\epsilon) |z_1|^2 + |z_4|^2 = |z_2|^2 + |z_3|^2, $$
and the inverse image of the moment map over $\XX_{\lambda_1}$ is the complex proto-manifold (diffeology) associated to two copies of the Kronecker foliation, namely the complex non-commutative torus $\Torus^2_\theta$ for $\theta = 1+ \epsilon$.

\subsection{} Let $P$ be a simple convex polytope. A $d$-dimensional convex simplicial polytope $Q$ is one where every proper face is a simplex. The dual $Q=P^\vee$ of a convex simple polytope $P$ is simplicial.  The face vector ($f$-vector) of $Q$ is defined to be $f(Q) = (f_{-1}=1, f_0,\ldots, f_{d-1})$ where $f_i$ is the number of $i$-dimensional faces of $Q$. It is an interesting question to characterize combinatorially all the vectors $f \in \integer^d$ that appear this way. 


We define the $h$-vector $h(Q)$ of  $Q$ by means of the following formula: 
$$h_i:=\sum_{j=0}^i  { {d-j} \choose {d-i} } (-1)^{i-j} f_{j-1}.$$

The $h$-vector of a simple polytope satisfies the  Dehn-Sommerville relations: $$h_i =h_{d-i},$$ that makes one think of Poincar\'e duality.

We further define the $g$-vector $g(Q)=(g_0,\ldots, g_{\lfloor d/2\rfloor})$ by means of $g_0:=h_0, g_1:=h_1-h_0, g_2 = h_2-h_1, \ldots, g_{\lfloor d/2\rfloor} :=h_{\lfloor d/2\rfloor} - h_{\lfloor d/2\rfloor -1}. $

The McMullen $g$-conjecture (1971) is a  characterization of the possible $f$-vectors for $Q=P^\vee$ in terms of the $g$-vectors. Shortly we will need the concept of an $M$-vector (cf. \cite{Stanleycommutative} page 55).  Given two integers $l, i>0$ there is a unique expansion of the form:
$$l = {{n_i} \choose {i}} + {{n_{i-1}}\choose {i-1}} + \cdots + {{n_j}\choose{j}},$$
where $n_i>n_{i-1}>\cdots > n_j \ge j \ge 1$. We define $l^{(i)}$ as:
$$l^{(i)}:={{n_i} \choose {i+1}} + {{n_{i-1}}\choose {i}} + \cdots + {{n_j}\choose{j+1}}.$$ A vector $(l_0,l_1,\ldots)$ is said to be an \emph{$M$-vector} if and only if $l_0=1$ and $0\le l_{i+1} \le l_i^{(i)}, i\ge1.$

The McMullen $g$-conjecture \cite{McMullenG} (1971) is a characterization of the possible $f$-vectors for $Q=P^\vee$ in terms of the $g$-vectors:

\begin{theorem} Let $h=(h_0,\ldots h_d) \in \integer^d$. The following two conditions are equivalent:
\begin{itemize}
\item There exists a simplicial $d$-polytope $Q$ such that $h(Q)=h.$
\item $h_0=1$, $h_i=h_{d-i}$ for all $i$, and $g(Q)$ is an $M$-vector.
\end{itemize}
\end{theorem}

The sufficiency of the McMullen conjecture was proved by Billera and Lee in \cite{BilleraLee}.

Famously, Richard Stanley (1983) proved that the necessity of the McMullen conjecture would be proved if one showed that, given $P$, there is a graded commutative algebra $R=R_0\oplus R_1\oplus \cdots$ over $\complex$ with $R_0=\complex$ generated by $R_1$ and with $\dim R_i = h_i - h_{i-1}$ for all $1\leq i \leq [d/2].$

He observed that this would be true in view of the hard Lefschetz theorem if there was a variety $X_P$ so that $\beta_{2i}(X_P) = h_i$, and $R=H^*(X_P)/(\omega)$ where $\omega$ is the hyperplane section (cf. page 76 of \cite{Stanleycommutative}).

The theory of toric varieties provides such an $X_P$ when $P$ is a rational simple polytope (Delzant polytope), and so Stanley proved the necessity of the McMullen conjecture for rational simplicial polytopes $Q$.

In the general case, the theory of non-commutative toric varieties provides us with (an albeit non-commutative) toric variety $\XX_P$ for an arbitrary (possibly irrational) simple polytope $P$. One must verify first that non-commutative toric varieties satisfy the hard Lefchetz theorem, a fact that we prove in \cite{NCTorics}: We have found a diffeological version of Chern's proof of the Hard Lefschetz theorem, but compare with \cite{battaglia2011foliations}. Hence we can give yet another proof of the (necessity of the) McMullen conjecture in the general case of an arbitrary simplicial polytope $Q$.

\subsection{} Let us consider the case $n-2m-1=2$ in which $P$ is a 2-dimensional polygon, and $\XX_P$ is of complex dimension $2$. In the rational case, it is known that, considering the unique (finite) Hirzebruch-Jung continued fraction expansion: 
$$a_1 -\cfrac{1}{a_2-\cfrac{1}{a_3-\cfrac{1}{\cdots-\cfrac{1}{a_r}}}}$$
of the slopes of the generating rays gives a procedure to subdivide the cones of the fan (or equivalently cutting off corners of the polygon) increasing at each step the number of sides by one and doing a blow-up at the level of the $X_P$, until all singularities are resolved.

In the case of a non-commutative toric variety $\XX_P$, some of the slopes will be irrational, and the continued fraction expansion will still be unique but infinite: 
$$a_1 -\cfrac{1}{a_2-\cfrac{1}{a_3-\cfrac{1}{a_4-\cfrac{1}{\cdots}}}}$$
Therefore, we can always push the non-commutative part of $\XX_P$ to smaller and smaller cones until, eventually, we get sort of a commutative resolution: The non-algebraic (but still complex) ``toric variety" $X_{P'}$ associated to a polygon $P'$ with an infinite (but countable) number of sides.

\section{Relation to Mirror Symmetry}

\subsection{} Over the last several years, a variety of new categorical structures
have been discovered by physicists.  Furthermore, it has become
transparently evident that the categorical language is 
suited to describing cornerstone concepts in modern theoretical
physics. Two major physical sources for categorical development are
the subject of Mirror Symmetry and the study of Topological States of
Matter.

\subsection{Mirror symmetry} There is a physical duality between $N= 2$
superconformal field theories.  In the 90's, Maxim Kontsevich
re-interpreted this concept from physics as a deep and ubiquitous
mathematical duality now known as Homological Mirror Symmetry (HMS).
His 1994 lecture created a frenzy of activity in the mathematical
community which lead to a remarkable synergy of diverse mathematical
disciplines: symplectic geometry, algebraic geometry, and category
theory. In essence, HMS is a correspondence between algebro-geometric
and symplectic data.

\subsection{ Topological states of matter} This is a ground-breaking field in
modern physics. These new states of matter appear not as a result of
symmetry breaking as proposed by Landau, but instead for topological
reasons (although symmetries such as charge conservation, parity and
time-reversal also play an important role). Topological insulators
were theoretically predicted by Bernevig and Zhang and Fu and Kane and
have been observed experimentally. Later Kitaev proposed a more
general classification scheme based on K-theory which applies to both
topological superconductors and topological insulators, but only in
the case of non-interacting fermions. One application was the prediction
\cite{LK}  of
the existence of phantom and quasi phantom categories - geometric
nontrivial categories with trivial homology or Grothendieck group
respectively. This prediction was recently proved by B\"ohning, Graf
von Bothmer, Katzarkov and Sosna \cite{BBKPS}, B\"ohning, Graf von
Bothmer and Sosna \cite{BBS}, Orlov and Alexeev \cite{AO2012}, Gorchinskiy
and Orlov \cite{GO}, and Galkin and Shinder \cite{GS}.

According to the HMS conjecture, the  mirrors of manifolds $X$ are Landau--Ginzburg models, namely proper maps $W : Y \to \complex$, 
where $Y$ is a symplectic manifold.  
The category $\Db(X)$ is expected to be equivalent to
 the Fukaya--Seidel category of $ W:  Y  \to  \complex$, see e.g. \cite{Auroux2008}. 
Toric varieties is the first example where HMS was proven  
see e.g. \cite{Auroux2008} and \cite{Abouzaid2006}. In \cite{Auroux2008}
HMS was proved for  non-commutative deformations of the complex projective plane. 
The examples of non-commutative toric varieties constructed in this paper 
provide an ample amount of interesting examples where HMS can be tested.

Recently, new  techniques have been developed (cf. \cite{BDFKK}). These techniques are based on the correspondence between
deformations on the A side and VGIT quotients on the B side.
The examples described in this paper can be treated via these techniques and we intend to do so in a future paper \cite{NCTorics}. Here we pose some questions:

A recent work of Witten \cite{WIT} suggests that Fukaya - Seidel  category  
(following HMS $D^b$ on the B side) 
of the A side can be related to a WKB approximation applied to the solutions  of some differential equation. 
  
{\bf Question 1.} Find the differential equation describing via WKP approximation  Fukaya - Seidel  category of the mirrors of non-commutative toric varieties.

{\bf Question 1.} Do  Fukaya - Seidel  categories  of the mirrors of noncommutative toric varieties contain phantom categories in their semiorthogonal decompositions?

\bibliographystyle{amsplain}
\bibliography{NCToricsbib}

\begin{thebibliography}{10}

\bibitem{Abouzaid2006}
Mohammed Abouzaid, \emph{Homogeneous coordinate rings and mirror symmetry for
  toric varieties}, Geom. Topol. \textbf{10} (2006), 1097--1157 (electronic).
  \MR{2240909 (2007h:14052)}

\bibitem{AO2012}
Valery Alexeev and Dmitri Orlov, \emph{Derived categories of burniat surfaces
  and exceptional collections}, Mathematische Annalen (2012), 1--17.

\bibitem{audin2004torus}
Michele Audin, \emph{Torus actions on symplectic manifolds}, vol.~93,
  Birkhauser, 2004.

\bibitem{Auroux2008}
Denis Auroux, Ludmil Katzarkov, and Dmitri Orlov, \emph{Mirror symmetry for
  weighted projective planes and their noncommutative deformations}, Ann. of
  Math. (2) \textbf{167} (2008), no.~3, 867--943. \MR{2415388 (2009f:53142)}

\bibitem{BBCG2010}
A.~Bahri, M.~Bendersky, F.~R. Cohen, and S.~Gitler, \emph{The polyhedral
  product functor: a method of decomposition for moment-angle complexes,
  arrangements and related spaces}, Adv. Math. \textbf{225} (2010), no.~3,
  1634--1668. \MR{2673742 (2012b:13053)}

\bibitem{BDFKK}
Matthew Ballard, Colin Diemer, David Favero, Ludmil Katzarkov, and Gabriel
  Kerr, \emph{The mori program and non-fano toric homological mirror symmetry},
  arXiv preprint arXiv:1302.0803 (2013).

\bibitem{battaglia2011foliations}
Fiammetta Battaglia and Dan Zaffran, \emph{Foliations modelling nonrational
  simplicial toric varieties}, arXiv preprint arXiv:1108.1637 (2011).

\bibitem{BilleraLee}
Louis~J. Billera and Carl~W. Lee, \emph{A proof of the sufficiency of
  {M}c{M}ullen's conditions for {$f$}-vectors of simplicial convex polytopes},
  J. Combin. Theory Ser. A \textbf{31} (1981), no.~3, 237--255. \MR{635368
  (82m:52006)}

\bibitem{BBKPS}
Christian B{\"o}hning, Hans-Christian~Graf von Bothmer, Ludmil Katzarkov, and
  Pawel Sosna, \emph{Determinantal barlow surfaces and phantom categories},
  arXiv preprint arXiv:1210.0343 (2012).

\bibitem{BBS}
Christian B{\"o}hning, Hans-Christian~Graf von Bothmer, and Pawel Sosna,
  \emph{On the derived category of the classical godeaux surface}, Advances in
  Mathematics \textbf{243} (2013), 203--231.

\bibitem{BP2002}
Victor~M. Buchstaber and Taras~E. Panov, \emph{Torus actions and their
  applications in topology and combinatorics}, University Lecture Series,
  vol.~24, American Mathematical Society, Providence, RI, 2002. \MR{1897064
  (2003e:57039)}

\bibitem{connesnoncommutative}
Alain Connes, \emph{Noncommutative geometry}, Academic Press Inc., San Diego,
  CA, 1994. \MR{1303779 (95j:46063)}

\bibitem{CoxToric2011}
David~A. Cox, John~B. Little, and Henry~K. Schenck, \emph{Toric varieties},
  Graduate Studies in Mathematics, vol. 124, American Mathematical Society,
  Providence, RI, 2011. \MR{2810322 (2012g:14094)}

\bibitem{CrainicMoerdijk}
M.~Crainic and I.~Moerdijk, \emph{A homology theory for \'etale groupoids}, J.
  Reine Angew. Math. \textbf{521} (2000), 25--46. \MR{2001f:58039}

\bibitem{crainic2004cech}
Marius Crainic and Ieke Moerdijk, \emph{Cech-de rham theory for leaf spaces of
  foliations}, Mathematische Annalen \textbf{328} (2004), no.~1, 59--85.

\bibitem{de1997new}
Santiago~L{\'o}pez de~Medrano and Alberto Verjovsky, \emph{A new family of
  complex, compact, non-symplectic manifolds}, Bulletin of the Brazilian
  Mathematical Society \textbf{28} (1997), no.~2, 253--269.

\bibitem{Tsygan}
B.~L. Fe{\u\i}gin and B.~L. Tsygan, \emph{Additive {$K$}-theory and crystalline
  cohomology}, Funktsional. Anal. i Prilozhen. \textbf{19} (1985), no.~2,
  52--62, 96. \MR{MR800920 (88e:18008)}

\bibitem{fulton1993introduction}
William Fulton, \emph{Introduction to toric varieties.(am-131)}, vol. 131,
  Princeton University Press, 1993.

\bibitem{GS}
Sergey Galkin and Evgeny Shinder, \emph{Exceptional collections of line bundles
  on the beauville surface}, arXiv preprint arXiv:1210.3339 (2012).

\bibitem{NearlyBook}
A.~Gonzalez, E.~Lupercio, C.~Segovia, and B.~Uribe, \emph{Topological quantum
  field theories of simension 2 on orbifolds}, To appear.

\bibitem{GO}
Sergey Gorchinskiy and Dmitri Orlov, \emph{Geometric phantom categories},
  Publications math{\'e}matiques de l'IH{\'E}S (2013), 1--21.

\bibitem{hamilton2007quantization}
Mark~D Hamilton, \emph{The quantization of a toric manifold is given by the
  integer lattice points in the moment polytope}, arXiv preprint
  arXiv:0708.2710 (2007).

\bibitem{hector2011rham}
Gilbert Hector, Enrique Mac{\'\i}as-Virg{\'o}s, and Esperanza
  Sanmart{\'\i}n-Carb{\'o}n, \emph{De rham cohomology of diffeological spaces
  and foliations}, Indagationes Mathematicae \textbf{21} (2011), no.~3,
  212--220.

\bibitem{HectorLie}
Gilbert Hector, Enrique Mac{\'{\i}}as-Virg{\'o}s, and Antonio Sotelo-Armesto,
  \emph{The diffeomorphism group of a {L}ie foliation}, Ann. Inst. Fourier
  (Grenoble) \textbf{61} (2011), no.~1, 365--378. \MR{2828134 (2012e:57049)}

\bibitem{HKR}
G.~Hochschild, Bertram Kostant, and Alex Rosenberg, \emph{{Differential forms
  on regular affine algebras.}}, Trans. Am. Math. Soc. \textbf{102} (1962),
  383--408 (English).

\bibitem{iglesias2010orbifolds}
Patrick Iglesias, Yael Karshon, and Moshe Zadka, \emph{Orbifolds as
  diffeologies}, Transaction of the American Mathematical Society \textbf{362}
  (2010), no.~6, 2811--2831.

\bibitem{iglesias2008diffeology}
Patrick Iglesias-Zemmour, \emph{Diffeology}.

\bibitem{Karshon2008}
Yael Karshon and Masrour Zoghi, \emph{An effective orbifold groupoid is
  determined up to morita equivalence by its underlying diffeological
  orbifold}, Preprint (2008).

\bibitem{LK}
L.~Katzarkov, \emph{Homological mirror symmetry and algebraic cycles},
  Homological mirror symmetry, Lecture Notes in Phys., vol. 757, Springer,
  Berlin, 2009, pp.~125--152. \MR{2596637 (2011b:14082)}

\bibitem{NCTorics}
L.~Katzarkov, E.~Lupercio, L.~Meersseman, and A.~Verjovsky,
  \emph{Non-commutative toric varieties}, To appear.

\bibitem{kontsevich2008xi}
Maxim Kontsevich, \emph{{XI} {S}olomon {L}efschetz {M}emorial {L}ecture series:
  {H}odge structures in noncommutative geometry. {N}otes by {E}rnesto
  {L}upercio}, Contemp. Math \textbf{462} (2008), 1--21.

\bibitem{McMullenG}
P.~McMullen, \emph{The numbers of faces of simplicial polytopes}, Israel J.
  Math. \textbf{9} (1971), 559--570. \MR{0278183 (43 \#3914)}

\bibitem{meersseman2004holomorphic}
Laurent Meersseman and Alberto Verjovsky, \emph{Holomorphic principal bundles
  over projective toric varieties}, Journal fur die Reine und Angewandte
  Mathematik (2004), 57--96.

\bibitem{MoerdijkFoliation}
Ieke Moerdijk, \emph{Models for the leaf space of a foliation}, European
  {C}ongress of {M}athematics, {V}ol. {I} ({B}arcelona, 2000), Progr. Math.,
  vol. 201, Birkh\"auser, Basel, 2001, pp.~481--489. \MR{1905337 (2003g:57045)}

\bibitem{MoerdijkSurveyOrbifolds}
\bysame, \emph{Orbifolds as groupoids: an introduction}, Orbifolds in
  mathematics and physics (Madison, WI, 2001), Contemp. Math., vol. 310, Amer.
  Math. Soc., Providence, RI, 2002, pp.~205--222. \MR{2004c:22003}

\bibitem{Rieffel}
Marc~A. Rieffel, \emph{{$C\sp{\ast} $}-algebras associated with irrational
  rotations}, Pacific J. Math. \textbf{93} (1981), no.~2, 415--429.
  \MR{MR623572 (83b:46087)}

\bibitem{souriau1980groupes}
J~Souriau, \emph{Groupes diff{\'e}rentiels}, Differential geometrical methods
  in mathematical physics (1980), 91--128.

\bibitem{Stanleycommutative}
Richard~P. Stanley, \emph{Combinatorics and commutative algebra}, second ed.,
  Progress in Mathematics, vol.~41, Birkh\"auser Boston Inc., Boston, MA, 1996.
  \MR{1453579 (98h:05001)}

\bibitem{WIT}
Edward Witten, \emph{A new look at the path integral of quantum mechanics},
  arXiv preprint arXiv:1009.6032 (2010).

\end{thebibliography}

\end{document}